\documentclass[12pt]{article}

\usepackage{amsmath,amssymb,amsfonts} 

\usepackage{natbib} 
\usepackage{abstract} 
\usepackage{authblk} 
\usepackage[utf8]{inputenc} 

\usepackage{geometry}
\usepackage{graphicx} 
\usepackage[labelfont=bf]{caption}
\usepackage{subcaption}
\usepackage{algorithm}
\usepackage{algpseudocode}
\usepackage{hyperref}  
\usepackage{multirow} 
\usepackage{multicol} 
\usepackage{hhline} 
\usepackage{lscape}

\usepackage{tikz}
\usepackage{capt-of}
\usepackage{booktabs}

\usetikzlibrary{decorations.pathreplacing,calc}
\newcommand{\tikzmark}[1]{\tikz[overlay,remember picture] \node (#1) {};}
\newcommand*{\AddNote}[4]{%
    \begin{tikzpicture}[overlay, remember picture]
        \draw [decoration={brace,amplitude=0.5em},decorate,ultra thick,black]
            ($(#3)!(#1.north)!($(#3)-(0,1)$)$) --  
            ($(#3)!(#2.south)!($(#3)-(0,1)$)$)
                node [align=center, text width=5cm, pos=0.5, anchor=west] {#4};
    \end{tikzpicture}
}%

\DeclareMathOperator*{\argmin}{arg\,min}
\algdef{SE}[DOWHILE]{Do}{doWhile}{\algorithmicdo}[1]{\algorithmicwhile\ #1}%

\begin{document}

\title{\textbf{Learning to repeatedly solve routing problems}}
\author[1,2,3]{Mouad Morabit}
\author[1,2]{Guy Desaulniers}
\author[3,4]{Andrea Lodi}

\affil[1]{MAGI, Polytechnique Montr\'eal, Canada}%
\affil[2]{GERAD, Montr\'eal, Canada}
\affil[3]{CERC, Polytechnique Montr\'eal, Canada}
\affil[4]{Jacobs Technion-Cornell Institute\\
Cornell Tech, USA}
\affil[ ]{\textit {\{mouad.morabit, guy.desaulniers\}@gerad.ca}}
\affil[ ]{\textit {andrea.lodi@cornell.edu}}
\date{\vspace{-5ex}}
\maketitle
\begin{abstract}
In the last years, there has been a great interest in machine-learning-based heuristics for solving NP-hard combinatorial optimization problems. The developed methods have shown potential on many optimization problems. In this paper, we present a learned heuristic for the reoptimization of a problem after a minor change in its data. We focus on the case of the capacited vehicle routing problem with static clients (i.e., same client locations) and changed demands. Given the edges of an original solution, the goal is to predict and fix the ones that have a high chance of remaining in an optimal solution after a change of client demands. This partial prediction of the solution reduces the complexity of the problem and speeds up its resolution, while yielding a good quality solution. The proposed approach resulted in solutions with an optimality gap ranging from 0\% to 1.7\% on different benchmark instances within a reasonable computing time. 
\end{abstract}

\section{\label{sec:intro}Introduction}

In the recent years, the idea of integrating machine learning (ML) and combinatorial optimization (CO) has been greatly explored. CO problems are NP-hard and there are generally two approaches to solve them. The \textit{exact} methods that are guaranteed to obtain an optimal solution but can be computationally very expensive for large instances, and the \textit{heuristic} methods that trade off the optimality of the solution for a reasonable computing time. The idea of leveraging ML for the development of new heuristics has shown potential in many CO problems such as traveling salesman problem (TSP), capacited vehicle routing problem (CVRP), etc. It is true that most of these learned heuristics do not outperform highly optimized and specialized CO algorithms, especially for problems that have been extensively studied in the literature. Nevertheless, the ideas behind them provide a certain flexibility for adjustments and applications to other problems for which no good heuristic exists, or they can be integrated in already existing algorithms to speed them up.

In this paper, we focus on CO applications in which a problem is repeatedly solved, e.g., daily or hourly, or even within a shorter interval, by changing neither its structure (e.g., its constraints) nor even its size (i.e., its variables), but only the data that define each instance solved in the specific time interval. This is the case of applications in which the infrastructure whose operations must be optimized does not change, for example a fleet of vehicles that deliver goods or the power plants producing energy, but the demand of goods or energy changes. And it is also the case of real-time changes to the data due to disruptions in the infrastructure, for example arcs disappearing from a network (i.e., their capacity going to 0) due to accidents. 

Applications of this type might be difficult to solve in reasonable amount of time, or, more precisely, each instance in isolation might require a significant computational effort even if the solution method has been designed after intensively studying the characteristics of the CO problem. This is the theoretical consequence of NP-hardness, and, on the practical side, it is due to the fact that the solution methods are largely designed to be agnostic to the data. 

The goal of this paper is to propose a learned heuristic allowing a fast reoptimization of a CO problem after a slight modification of its data. In other words, we put ourselves in a way more restrictive context with respect to the use of ML for CO problems: we do not want to leverage ML to devise a heuristic that produces good feasible solutions for all, say, TSP instances within the same distribution. We are settling for a lesser objective, i.e., that of learning what can be (more or less with high probability) safely left unchanged in the solution of a reference instance of a CO problem when the data of the instance are perturbed. Indeed, the intuition is that given an instance and its solution $S_1$, if the instance is reoptimized after a slight change in its data, the new solution $S_2$ will have a significant overlap with $S_1$, while only some (minor) parts of the solutions will be different. Therefore, instead of reoptimizing from scratch, the goal is to predict the parts of the solution that have a high probability of remaining the same. The corresponding variables can be fixed, thus reducing the search space and accelerating the resolution of the problem, almost independently of the solution method applied. 

To provide a concrete example, let us consider the CVRP, which is the problem where we will apply the method proposed in this paper. The goal of the CVRP is to construct vehicle routes in order to serve geographically-dispersed clients while minimizing the travel costs and respecting vehicle capacity.
Let us take the example of a delivery company that solves CVRPs on a daily basis. For a given day, the company observes that the clients are the same as in the instance solved the previous day (i.e., same client locations) and that only some clients have a different demand. After the optimization of the problem, a similarity between the solutions is noticed. Given the optimal solution (or a heuristic one) already in hand and the new demands, the objective would be to predict and fix the sequences of edges that have a high probability to remain the same. In case of a graph-structured problem like the CVRP, it is also possible to reduce the size of the network by aggregating the nodes/edges of the fixed sequences, therefore accelerating the resolution of the problem and reducing its complexity. Note that the predictions obtained by the learned model are not necessarily 100\% accurate, misclassifications may occur and thus affect the quality of the solution. The goal is to find a good compromise between the quality of the solution (i.e., optimality gap) and the computing time.

As observed, one can think of several CO problems where slight modifications to the problem data lead to similar solutions after reoptimization, especially for problems that are solved repeatedly and for which a data set is already available (e.g., the unit commitment problem in which the power plants producing energy are always the same but the demand changes daily or hourly). The learned model will try to partially predict an optimal solution. In this paper, we will only consider the case of the CVRP with changing demands and fixed customers, but the method offers some flexibility and has the potential to be applied to other problems or to be integrated into existing algorithms.

The remainder of the paper is organized as follows. In Section 2, we present some recent work on using ML for solving CO problems or accelerating their solution process. Section 3 is devoted to presenting the CVRP, with a focus on the methods we consider to solve it. In Section 4, we cover all the details of the method we propose. Section 5 reports our computational results. Finally, conclusions are drawn in Section 6.

\section{\label{sec:related-work}Related work}

In the literature, several heuristics incorporating ML models for solving NP-hard CO problems were explored (see, the recent surveys \citet[]{bengio2018, fioretto, cappart2021combinatorial}). The proposed learning methods mostly fall into one of two categories. In the first category, \textit{supervised learning} methods (examples of the imitation learning paradigm) are algorithms that learn from data and try to mimic an expert. The data is given to the learner as a pair of features and expected outputs (or labels) and the learner tries to find patterns in the data while optimizing a performance measure. Generally, the aim of this approach is to replace known expensive computations by fast approximations (e.g., for our case, the expensive computation corresponds to reoptimizing the problem from scratch). In the second category, we find the \textit{reinforcement learning} (RL) algorithms that apply a ``learning-by-experience" paradigm. Instead of giving the learner the data on which to learn, RL algorithms explore the decision space by interacting with its environment in order to achieve a certain goal. In response to a decision (i.e., an action), the learner receives a real-valued reward. The goal is to find the best decisions to make at each state while maximizing the expected rewards. This learning approach has the advantage of not requiring any data and has generally shown a better generalization, i.e., the ability of continuing to be effective when the problems change for example in size or even better in data distribution.

Several studies have tried to tackle the TSP, the most classical CO problem where the goal is to visit a set of nodes exactly once with a single vehicle while minimizing the travel distance. For the supervised approaches, we mention the work of \citet{joshi2019}, where a learned heuristic method is presented. The method is based on a graph convolutional network model that takes the entire graph as an input and outputs an adjacency matrix with associated edge probabilities, which are then used to build a valid tour using a beam search algorithm. The authors report good results on fixed size instances. However, a very poor generalization is noticed when testing the models on instances of different size. Instead, the methods based on RL have shown more potential. The idea of using RL to solve CO problems was explored in \citet{neural-co} with an application to the TSP and the knapsack problem. Further studies were then conducted, such as \citet{attention-routing}, where the authors present an encoder-decoder model based on the attention mechanism \citep{attention-all-u-need}. In a similar line, \citet{reinforcement-for-vrp} present a framework focusing on solving the CVRP, which can be considered as a generalization of the TSP for the case with multiple vehicles. The results reported by the authors show a better performance when compared to the OR-Tools solver and other heuristic algorithms. Other researchers contributed to the methods for solving the broader class of graph-structured CO problems. \citet{dai-co-graphs} initiated the idea of learning on graphs, and the works of \citet{li-co-gcnn} and \citet{sahil-large-graph-drl} enhanced further the scalability to larger graphs. These methods have been applied to various NP-hard graph problems such as the minimum vertex cover, maximum clique, influence maximization problem, etc.

Unlike the previous works that seek to build an end-to-end solution to the different problems, other methods focus on using ML to guide and accelerate the solution process. Since we pay special attention to mixed-integer programs (MIP), it is worth mentioning the various projects involving ML in the context of branch and bound \citep[]{khalil2016,alvarez2017,gasse2019,zarpellon-param-bab}, many of which seek to learn a branching policy imitating the strong branching strategy. Another potential use of ML is embedding a learned model in MIP solvers in order to decide if a decomposition of the problem is beneficial or not \citep{kruber-decomp}, or if it is favorable to linearize the quadratic part of a mixed-integer quadratic program \citep{bonami-quad}.

Other learning-based methods have proven to be very effective on problems that are solved repeatedly, especially when the input data changes only slightly, which is the key idea of our project. These methods can exploit existing data from previous solutions in order to speed up the resolution of similar unseen instances. \citet{xavier-unit-commit} exploit the idea and apply it on the security-constrained unit commitment, a problem occuring in power systems and electricity markets. The authors report high speedups on computing time, up to 10 times faster than solving the problem from scratch, and without a noticeable loss in solution quality. Along the same lines, \citet{lodi-facility-location} consider an application to the facility location problem. They seek to estimate the proportion of a solution that has a high probability of remaining unchanged after a perturbation. An additional constraint is added to the original formulation according to the predictions obtained by a regression model. Both \citet{xavier-unit-commit} and \citet{lodi-facility-location} leverage ML for speeding up the resolution of repeatedly solved problems and do not seek to build an end-to-end solution, which is similar to our case. As opposed to \citet{xavier-unit-commit}, our method is related to that of \citet{lodi-facility-location} in the sense that both assume a reference solution to which changes (perturbations) are applied. However, the main difference lies on the fact that our approach revolves around fixing parts of the reference solution (i.e., edges) instead of estimating a bound on the number of changes without specifying which variables to set. Other than that, the CVRP remains a quite different problem and the approaches developed in this paper have the potential to be extended to other CVRP variants and also to other routing problems.

\section{\label{sec:cvrp}The capacited vehicle routing problem}

The CVRP is a CO problem that has been studied for many years, resulting in several exact and heuristic methods to solve it. It is classified as an NP-hard problem and remains a difficult problem to solve to optimality even with just a few hundred clients. Given a fleet of vehicles assigned to a depot, the problem consists in determining a set of possible routes (i.e., one route per vehicle used) to deliver goods to a set of dispersed clients while minimizing the travel costs. A route starts from the depot and visits a sequence of clients before returning back, and is considered feasible if the total amount of goods delivered does not exceed the vehicle capacity $Q$.

For solving this problem, we consider the algorithms based on column generation (CG - \citet{desaulniers2005a}), which is an exact iterative method for solving large linear programs. To ensure the obtention of integer solutions, CG is often embedded in a branch-and-bound framework where the linear relaxation of the problem is solved at each node using CG. In this case the method is referred to as \textit{branch-and-price} (B\&P) and it is considered the state-of-the-art exact method for solving the CVRP. But due to the complexity of the problem, solving large instances to optimality can be computationally expensive. Therefore, several heuristics have been proposed in the literature. 

Note that our approach does not necessarily require learning from optimal solutions. Since we want to apply it on instances of a reasonable size (i.e., 100 clients and more), solving the instances to optimality in order to collect data can be time consuming. Therefore, we chose to use a heuristic in the data collection phase, but during the evaluation phase, we exploit an exact B\&P algorithm. In the next sections, we present the problem formulation along with both the exact and the heuristic algorithms used.

\subsection{\label{sec:cvrp-formulation}CVRP formulation}

In this section, we formulate the CVRP as a set partitioning problem. Let $C$ be the set of clients to be serviced, $\Omega$ the set of all feasible routes and $c_r$ the cost of a route $r \in \Omega$. We define $a_{i}^r$ as a binary parameter equal to $1$ if client $i \in C$ is serviced by route $r \in \Omega$ and $0$ otherwise. Let $\theta_r$ be a binary decision variable equal to $1$ if route $r$ is part of the solution and $0$ otherwise. The problem can therefore be formulated as follows:

\begin{align}
(P) \hspace{5mm} \min_{\theta}\hspace*{3mm} \sum_{r\in \Omega} c_r\theta_r & \label{cvrp-obj} \\
\mbox{s.t.} \hspace*{3mm} \label{const-partition}
\sum_{r\in \Omega} a_{i}^r \theta_r = 1 &,\hspace*{3mm} \forall i \in C,  \\ \label{const-int}
\theta_r \in \{0,1\}&,\hspace*{3mm} \forall r \in \Omega,
\end{align}
where the objective~(\ref{cvrp-obj}) minimizes the total cost of the routes. Constraints~(\ref{const-partition}) ensure that each client is visited exactly once and constraints~(\ref{const-int}) are the binary requirements on the decision variables $\theta_r$. 

One can notice that for large instances, the size of the route set $|\Omega|$ becomes prohibitively large and it would not be possible to enumerate all the variables of the problem. This is why a CG-based algorithm is used. The goal is to start with a subset of variables and generate potentially improving columns when necessary.

\subsection{\label{sec:branch-and-price}Exact branch-and-price algorithm}

B\&P algorithms \citep{barnhart1998} are considered state-of-the-art exact algorithms for solving a variety of optimization problems (e.g., routing, scheduling, \dots ). B\&P is based on the branch-and-bound method in which the linear relaxation at each node is solved using CG. The CG process is iterative and consists in alternating between the resolution of a restricted version of the original linear relaxation, called a restricted master problem (RMP), and a pricing problem (PP). The goal of the PP is to find new improving columns of negative reduced cost (for a minimization problem) that can be added to the RMP. The CG process stops when no such columns are found. 
A branching then occurs and the branch-and-bound tree exploration continues. Optionally, cutting planes can also be added to strengthen the relaxation at each node, resulting on what is called a branch-cut-and-price method.

\subsubsection{\label{sec:bap-rmp}The restricted master problem}

The RMP corresponds to the linear relaxation of the formulation~(\ref{cvrp-obj})-(\ref{const-int}) but limited to only a subset $\mathcal{R} \subset \Omega$ of the variables. It is solved at each CG iteration, yielding a pair of primal and dual solutions $(\theta, \pi)$. The dual values $(\pi_i)_{i \in C}$ associated with the constraints~(\ref{const-partition}) are then used to find routes $r \in \Omega \backslash \mathcal{R}$ of negative reduced cost by solving the PP. If none exist, the CG process stops and the solution to the current RMP is thus optimal for the whole linear relaxation. Otherwise, the routes are added to the RMP (i.e., to the subset $\mathcal{R}$) which is then reoptimized.

\subsubsection{\label{sec:bap-pp}The pricing problem}

The PP can be defined as $min_{r \in \Omega}\thickspace\{ c_r - \sum_{i \in C} a_i^r \pi_i \}$. For many applications, especially routing and scheduling problems, this problem can be modeled as an elementary shortest path problem with resource constraints (ESPPRC - \citet{irnich-spprc}) where the goal is to find the least cost path between the source and destination nodes while visiting the nodes at most once (i.e., elementarity requirements) and respecting the resource constraints. For the CVRP, the only resource is the load, and a path is considered feasible if the load does not exceed the vehicle capacity. This problem can be defined over a graph $G=(V,A)$, where $V$ is the set of nodes representing the clients, in addition to the depot nodes $s$ and $t$, i.e., the source and destination nodes, respectively. The set $A$ represents the arcs, where each arc has an associated cost $c_{ij}, (i,j) \in A$. Hence, the cost of a path in $G$ is given by the sum of the costs $c_{ij}$ of its arcs. 

In order to take into account the dual values obtained by the RMP, at each iteration and for each arc $(i,j) \in A$, a modified cost $\bar{c}_{ij} = c_{ij} - \pi_i$ is used instead, where $\pi_i$ are the duals associated with constraints~(\ref{const-partition}) and $\pi_s = 0$. This guarantees that the cost of a feasible route in the network is equal to its reduced cost:
\begin{align}
\bar{c}_r = c_r - \sum_{i \in C} a_i^r \pi_i = \sum_{(i,j)\in A}c_{ij}b_{ij}^r - \sum_{i \in C} a_i^r \pi_i = \sum_{(i,j)\in A}(c_{ij} - \pi_i)b_{ij}^r = \sum_{(i,j) \in V} \bar{c}_{ij}b_{ij}^r
\end{align}
where $b_{ij}^r$ is equal to $1$ if arc $(i,j)\in A$ is traversed in route $r \in \Omega$, $0$ otherwise.

The ESPPRC is an NP-hard problem, which is mainly due to the client elementarity requirements since negative cost cycles can exist when using the modified costs $\bar{c}_{ij}$. Relaxing this constraint (i.e., allowing the generation of paths with cycles) leads to the SPPRC, which is an easier problem that can be solved in pseudo-polynomial time but yields a lower bound of inferior quality if solved as the PP. Other alternatives based on relaxations have been proposed in the literature, e.g., SPPRC-k-cyc \citep{irnich-spprc-k-cyc} and ng-routes \citep{baldacci-ng-routes-2011}. 

The success of B\&P methods for solving the CVRP and other variants is in good part due to the efficient methods for solving the PP, mainly using dynamic programming. A labeling algorithm is commonly used where a label corresponds to a partial path in $G$ (not to confound with the notion of label in ML). The algorithm starts with an initial label representing the trivial path containing only the source node $s$, it then gets extended forward along the outgoing arcs until reaching the destination node $t$. A new label is created at each extension if it yields a feasible path. At the end of the algorithm, the labels at the destination node $t$ representing negative reduced cost routes are used to build the new columns that are added to the RMP. Generally, several routes are added at once, which is known to speed up the solution process and to reduce the number of CG iterations. 

The CVRP remains one of the most studied CO problems in the literature. An efficient B\&P algorithm can be quite sophisticated and may contain several components. Discussing all the details is beyond the scope of this paper and we refer the interested reader to the survey of \citep{costa-contardo-desaulniers-2019} for an in-depth overview of the methods.

\subsubsection{\label{sec:bap-impl}Branch-and-price implementation}

In this work, we consider using \textit{VRPSolver} \citep{vrpsolver-pessoa}, which is a generic implementation of an exact branch-cut-and-price method for VRP problems. The advantage of using VRPSolver lies in the fact that it combines several algorithms and acceleration techniques introduced by several authors, e.g., ng-routes, path enumeration, bi-directional labeling, stabilization, etc. Implementing these techniques from scratch would be otherwise very time consuming. The authors report excellent performance on several benchmark instances of various CVRP variants (e.g., with time windows, heterogeneous fleet, multiple depots, pickups and deliveries, etc). For more information about the implementation and the algorithms included in VRPSolver, the reader is referred to \citet{vrpsolver-pessoa}.

\subsection{\label{sec:cvrp-heuristic}Heuristic algorithm}

As previously mentioned, the CVRP remains a difficult problem to solve to optimality and can be time consuming when working with large instances. For our method, we need to solve several instances to collect enough data for the training phase, but solving them to optimality is not required (although recommended), since it is possible to obtain very high-quality solutions using specialized CVRP heuristics. For this reason, we chose to use a recent heuristic called FILO \citep{filo}, which is a short term for \textit{Fast Iterated Localized Search Optimization}. 

The method is based on the iterated local search paradigm and is specifically designed to solve large-scale instances of the CVRP. The algorithm starts by constructing an initial feasible solution using an adaptation of the savings algorithm by \citet{clarke-wright}. It is followed by an optional step that aims at reducing the number of vehicles used in the initial solution, if the latter is larger than a computed estimate (using a bin-packing greedy algorithm). The algorithm then proceeds to the core optimization step, which is based on a sequence of ruin and recreate steps. The \textit{ruin} step removes a certain number of vertices by means of a random walk of a given length (the ruin intensity can be controlled by adapting the random walk length). Then, the \textit{recreate} step tries to reinsert the removed vertices while trying to improve upon the best solution found. This is achieved by means of a number of local search operators targetting the vertices involved in the disruptive effects of the ruin step. The algorithm tries to keep a good balance between intensification and diversification, i.e., the local search can be concentrated on the parts of the solution that have not seen any improvement after several attempts, and at the same time, a continuous diversification is intended in order to escape from local optima. The algorithm stops after a determined number of iterations and the best solution found is returned.

The performance of the FILO algorithm has been compared by the authors to other state-of-the-art heuristics and has proven to be highly competitive on the X instances introduced by \citet{uchoa-x-instances}, which are also the benchmark instances used by our method. Note that the purpose of this section is not to present or compare the different heuristics that exist for solving the CVRP, but simply to present the heuristic that we used to obtain high-quality solutions during the data collection phase (i.e., with optimality gaps of less than $0,1\%$ on average). The method has also the advantage of having an open-source implementation that is freely accessible.

\section{\label{sec:methodology}Methodology}

The goal of this project is to accelerate the reoptimization of repeatedly solved CO problems for which there is only a slight change in the problem data. Let $\mathcal{P}_o$ and $\mathcal{S}_o$ be a problem instance and a computed good-quality solution, respectively. Furthermore, let $\mathcal{P}_m$ be a modified instance obtained by applying minor changes to $\mathcal{P}_o$ and for which no solution is known. Instead of optimizing $\mathcal{P}_{m}$ from scratch, the objective is to identify the parts of $\mathcal{S}_o$ that have a high probability of also being part of a solution to $\mathcal{P}
_{m}$ denoted by $\mathcal{S}_{m}$. 

\begin{figure}
     \centering
\subfloat[\small Solution of $\mathcal{P}_o$\label{subfig:original_solution}]{\tikz[remember
picture]{\node(1A){\includegraphics[width=0.35\textwidth]{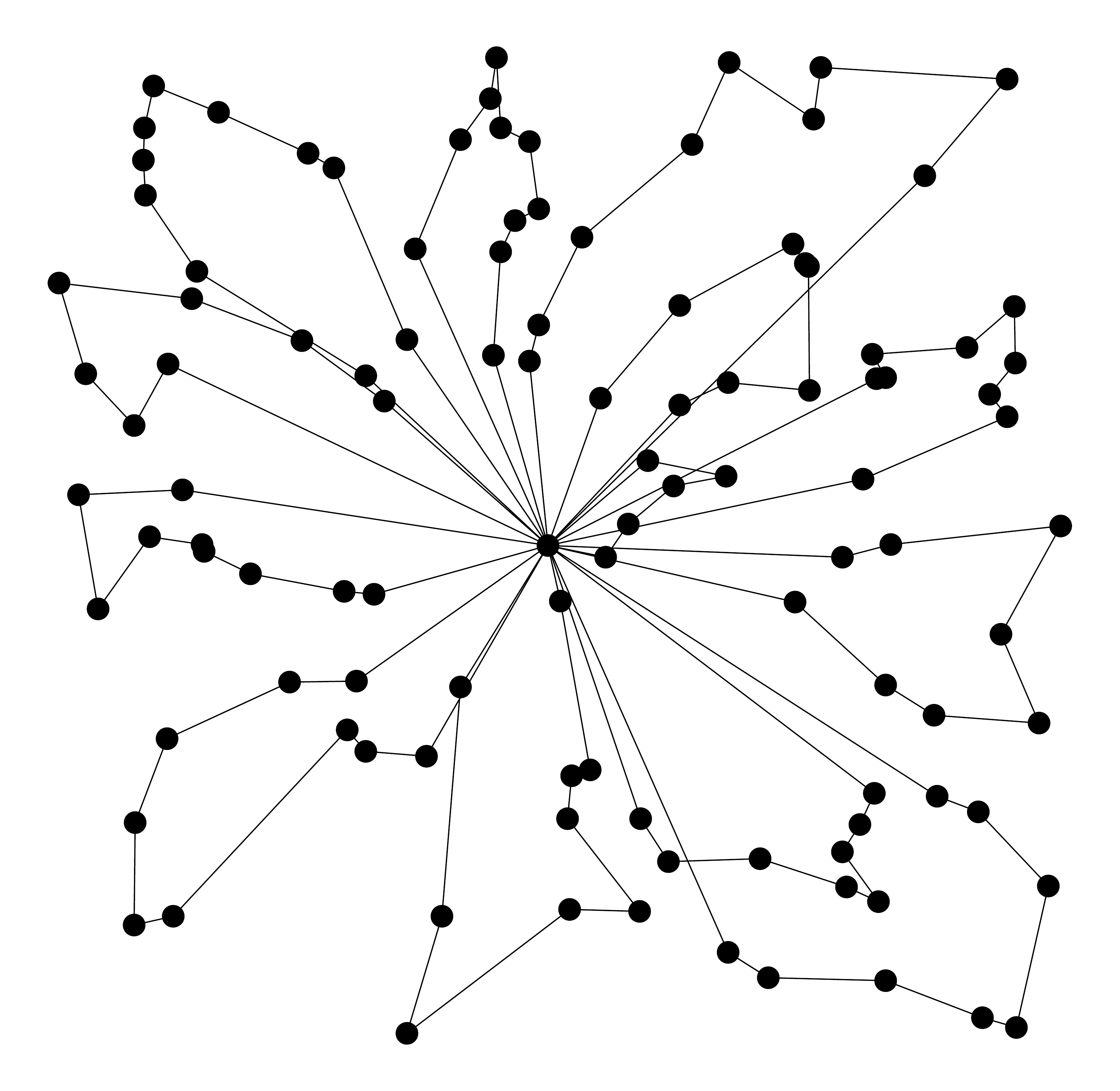}};}}%
\hspace*{0.25\textwidth}%
\subfloat[\small Solution of $\mathcal{P}_{m}$\label{subfig:modified_solution}]{\tikz[remember picture]{\node(1B){\includegraphics[width=0.35\textwidth]{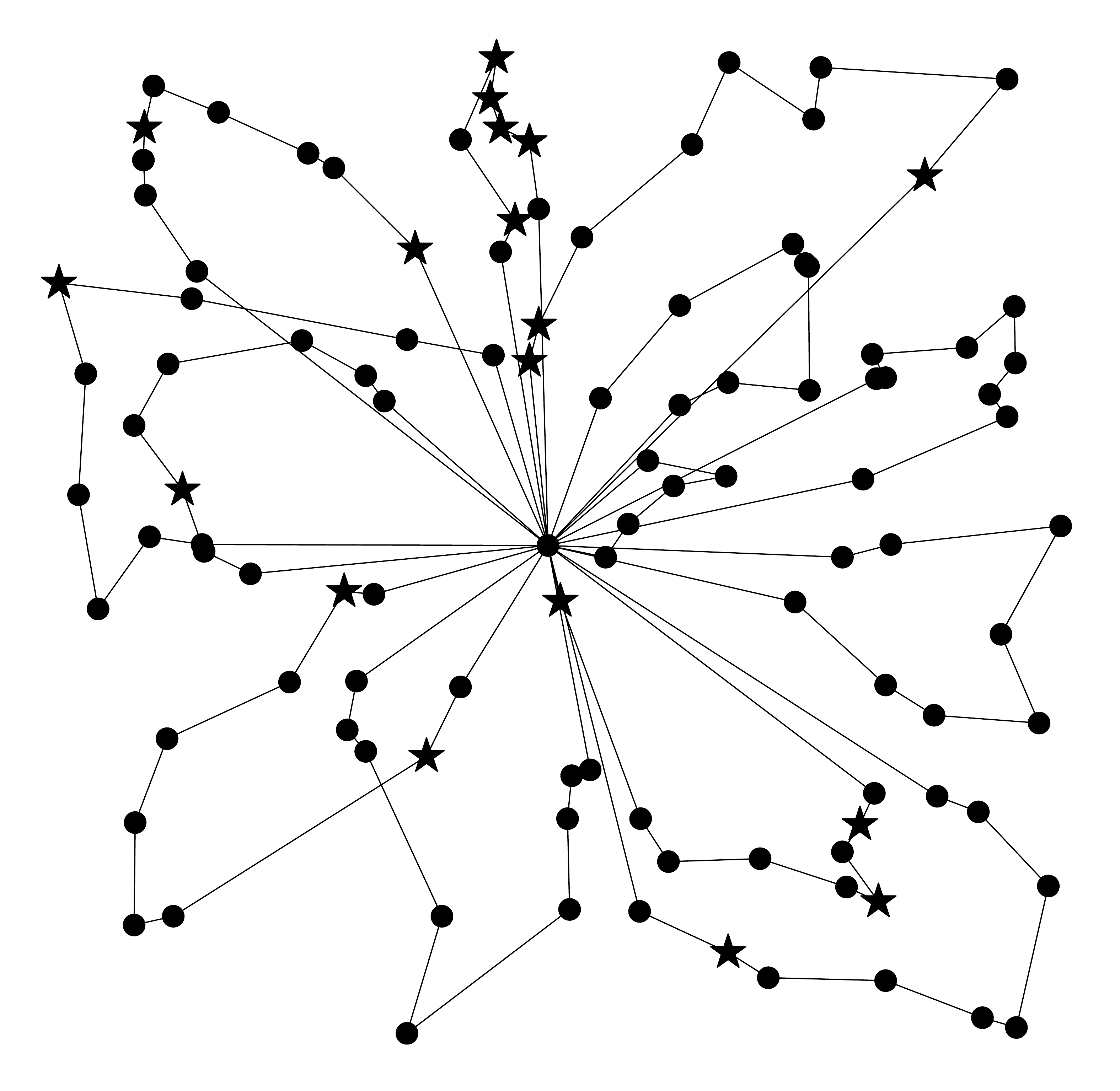}};}}

\tikz[overlay,remember picture]{\draw[-latex,thick] (1A) -- (1A-|1B.west)
node[midway,below,text width=0.2\textwidth]{After changing demands};}
\subfloat[\small Edges to predict and fix\label{subfig:similarity_solutions}]{\tikz[remember
picture]{\node(1C){\includegraphics[width=0.35\textwidth]{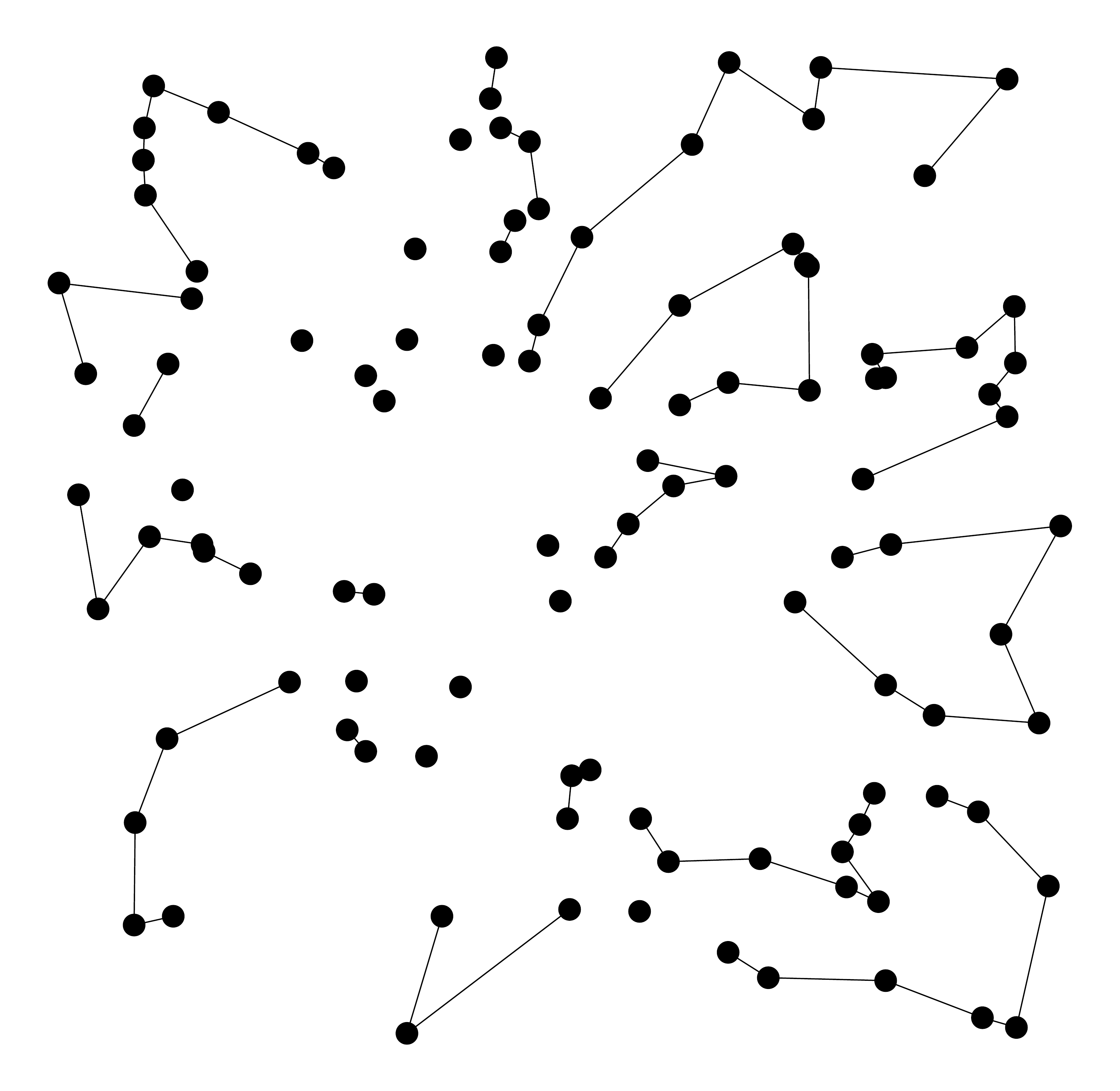}};}}%
\hspace*{0.25\textwidth}%
\subfloat[\small Edges to fix after aggregating sequences\label{subfig:aggregated_edges}]{\tikz[remember picture]{\node(1D){\includegraphics[width=0.35\textwidth]{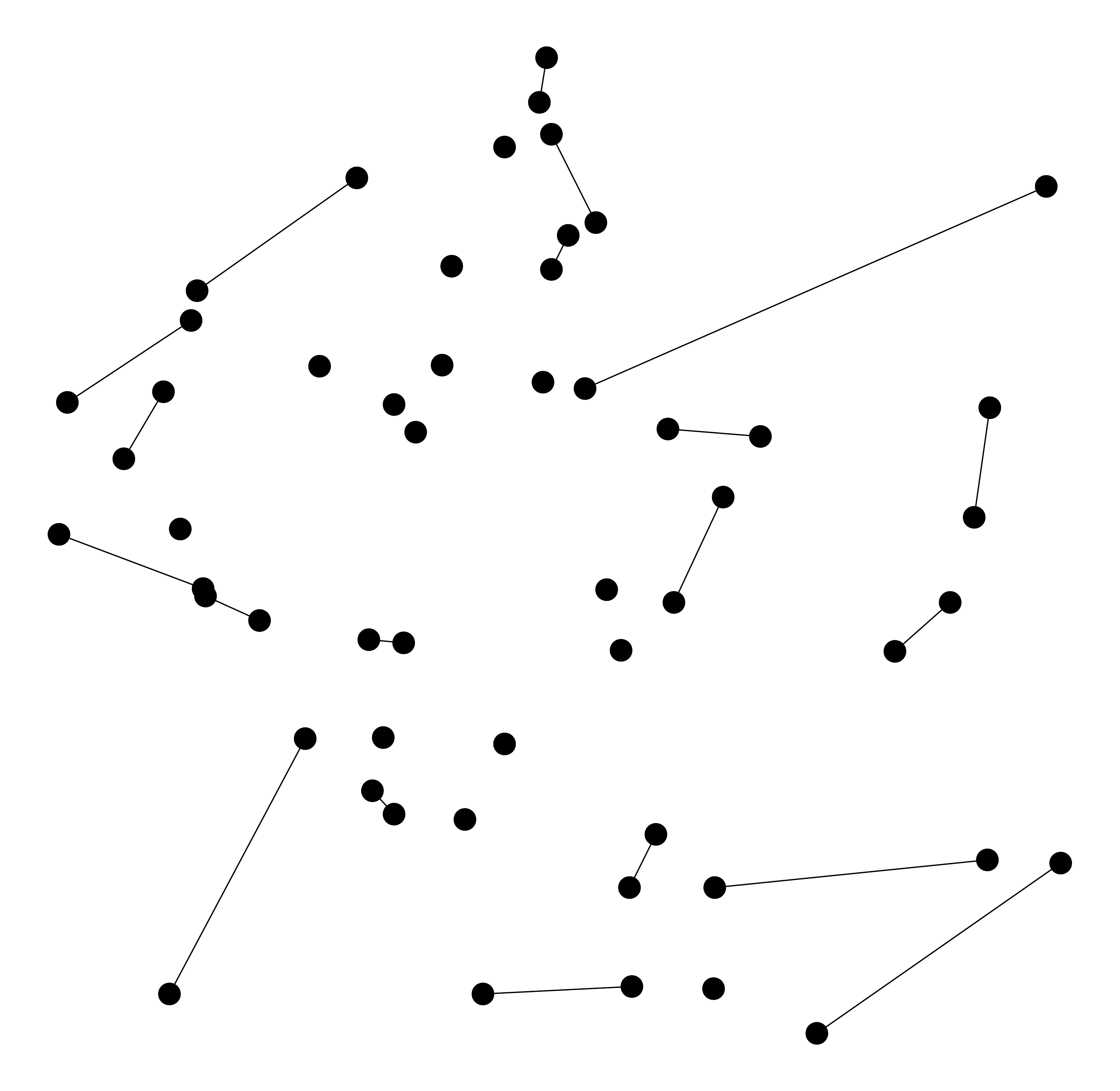}};}}

\tikz[overlay,remember picture]{\draw[-latex,thick] (1C) -- (1C-|1D.west)
node[midway,below,text width=0.2\textwidth]{After aggregating sequences};}
        \begin{center}
        \caption{Overview of the different steps of the method.}
        \label{fig:comparison_solutions_afterchange}
        \end{center}
\end{figure}

In this paper, we focus on the CVRP where the locations of the clients are the same but the demands are slightly different. Because we assume that the fleet of vehicles is unlimited and the travel costs are symmetric, we consider for the rest of this paper an undirected graph and denote by $E(\mathcal{S}_{o})$ the set of edges used in the known solution $\mathcal{S}_{o}$ to the original instance $\mathcal{P}_o$. Since the solution consists of vehicle routes, if we consider $E(\mathcal{S}_o)$ as the set of edges used in a solution $\mathcal{S}_o$ already in hand, the method aims at predicting the edges $e \in E(\mathcal{S}_o)$ that have a high chance to also be part of $\mathcal{S}_{m}$. By fixing these edges, a significant reduction of the problem complexity can be achieved, thus greatly reducing the computing time. Let us take the example illustrated in Figure \ref{fig:comparison_solutions_afterchange}, where we can observe the following: Figure (\ref{subfig:original_solution}) represents a solution $\mathcal{S}_o$ to an original CVRP instance involving $110$ nodes. The central node corresponds to the depot and the other nodes represent the clients. Figure (\ref{subfig:modified_solution}) shows a solution $\mathcal{S}_{m}$ to an instance $\mathcal{P}_{m}$ obtained by randomly changing the demands of $20\%$ of the clients of $\mathcal{P}_o$ (the clients with a changed demand are marked with a star node). By comparing the two figures, one can notice a significant similarity between the two solutions. Assuming we do not have yet the solution $\mathcal{S}_{m}$, if we succeed in predicting and fixing the edges $e \in E(\mathcal{S}_{o})$ that will be part of $\mathcal{S}_{m}$, a partial solution can be provided to the solver before starting the optimization. This is shown in Figure (\ref{subfig:similarity_solutions}), where the overlapping parts of the two solutions are highlighted, i.e., the edges $e \in E(\mathcal{S}_{o}) \cap E(\mathcal{S}_{m})$. Furthermore, it is also possible to reduce the size of the network by aggregating the sequences comprised of two or more edges into a single one as shown in Figure (\ref{subfig:aggregated_edges}). As a result, the number of nodes is reduced from $110$ to only $54$.

On the other hand, since the predictions obtained by the ML model are not always accurate and errors may occur, fixing the wrong edges can affect the quality of the solution obtained. In fact, fixing many edges implies shorter computing times but more chances to fix the wrong ones. On the other hand, if the number of fixed edges is small, there is a greater chance of obtaining a better quality solution though with a higher computing time. By controlling the number of fixed edges (e.g., by tuning the hyperparameters of the model), it is possible to find the right compromise between the quality of the solution and the computing time.

\subsection{\label{sec:data-collect}Data collection}

We chose to address this learning problem using a supervised learning approach. More precisely, a binary classification model is employed. The first step in the process is to collect enough data for the training. Given a tuple of original and modified instances and their solutions $(\mathcal{P}_{o},\mathcal{S}_{o},\mathcal{P}_{m},\mathcal{S}_{m})$,
a labeled dataset $\mathcal{D} = \{\{\boldsymbol{x}_e, y_e\} | \forall e \in E(\mathcal{S}_{o})\}$ is built where each entry represents an edge in the original solution, the vector $\boldsymbol{x}_{e} \in \mathbb{R}^n$ corresponds to the edge features (i.e., input), where $n$ is the number of features and $y_e = \{0,1\}$ is the desired output (i.e., label). One can notice that we are only interested in the edges of $\mathcal{S}_{o}$ and not all the edges of the graph. 

\textbf{The labels.} Since we tackle this problem in a supervised manner, we need both solutions $\mathcal{S}_{o}$ and $\mathcal{S}_{m}$ to build the dataset, i.e., we need to give both the input and the desired output to the learner. 
The labels are assigned by simply checking the overlapping edges between the solutions $\mathcal{S}_{o}$ and $\mathcal{S}_{m}$ as follows:
\begin{align}
y_e = 
    \begin{cases}
       1 & \text{if $e \in E(\mathcal{S}_{o}) \cap E(\mathcal{S}_{m}) $}\\
      0 & \text{otherwise}
    \end{cases}
    , \hspace{5mm} e \in E(\mathcal{S}_{o}).
\end{align} \vspace{1mm}

\textbf{The features.} The features $\boldsymbol{x}_{e}$ represent the characteristics of each edge $e = \langle i,j \rangle$ in the original solution. The extracted features are the following:
\begin{itemize}
    \item The $(x,y)$ coordinates of both nodes $i$ and $j$;
    \item The cost of edge $c_e$;
    \item The old and new demands of nodes $i$ and $j$ (the demand is set to $0$ for the depot nodes $s$ and $t$);
    \item The distance between the depot and nodes $i$ and $j$;
    \item A boolean value indicating whether the edge is a depot edge, i.e., equal to $1$ if either $i$ or $j$ is a depot node;
    \item A boolean value indicating whether the client $i$ or $j$, or both, have a changed demand;
    \item The rank of $i$ with respect to $j$ (and vice versa) according to the neighbor distances, e.g., if $j$ is the nearest neighbor to $i$, then its rank is $1$, if $j$ is the second closest neighbor to $i$ then its rank is $2$, and so on. 
\end{itemize}

Note that since we are working on a symmetric version of the CVRP and that the order of the entries in the vector $\boldsymbol{x}_{e}$ is important, we always assume that $i < j$.

\subsection{\label{sec:ml_heuristic}ML prediction}
 
Once the dataset is in hand, we apply common ML practices such as  data preprocessing (i.e., normalization, encoding) and data splitting (i.e., dividing the data into a training, a validation and a test set), etc. As mentioned before, the prediction task corresponds to a binary classification problem and different ML models can be considered, e.g., Random Forest, SVM, Neural network.
In this section, we assume that we already have a trained predictive model that takes as input the edge features $\boldsymbol{x}_e$ and outputs the predictions $\hat{y}_{e}$ and we wish to know what to do with those predictions afterwards (the comparison of which ML model is best performing is discussed in the next section).

\subsubsection{\label{sec:edge_fixing}ML-based edge fixing}
Let $G = (V,E)$ be an undirected graph, where $V$ is the set of nodes including the depot nodes and $E$ the set of edges. Let $c_{e}, e \in E$ be the cost of edge $e=\langle i,j \rangle$. The edge flow can be written in terms of the master problem (\ref{cvrp-obj})-(\ref{const-int}) variables as
\begin{align}
x_e = \sum_{r \in \Omega} b_{e}^{r} \theta_{r},
\end{align}
where $b_{e}^{r} \in \{0,1\}$ indicates whether edge $e \in E$ is part of a route $r \in \Omega$.
Therefore, an edge can be fixed by adding a new constraint to the master problem. For each edge in the original solution $\mathcal{S}_o$, the flow is set to $1$ depending on the predictions obtained by the ML model, namely
\begin{align}
x_e = 1 , \hspace{5mm} \forall e \in E(\mathcal{S}_o) \hspace{1mm}:\hspace{1mm} \hat{y}_e = 1.
\end{align}

\subsubsection{\label{sec:infeasibility_case}Infeasibility case}

Sometimes fixing edges can lead to an infeasible restricted problem. This may only occur in the case when the total demand of a sequence of clients fixed by the ML model exceeds the vehicle capacity $Q$. For each edge, in addition to the output $\hat{y}_e$ returned by the model, it is possible to obtain its probability estimate $\hat{p}_e$ of being fixed. Generally, the model assigns the class $1$ if the probability estimate $\hat{p}_e$ is greater than $0.5$ and $0$ otherwise. For a sequence of clients whose sum of the demands exceed $Q$, a possible solution would be to identify and unfix the edges with lowest probabilities, until obtaining one or several feasible subsequences.

Let $G_{\hat{y}}=(V_{\hat{y}}, E_{\hat{y}})$ be the graph obtained by keeping only the fixed edges and their corresponding nodes, i.e., $E_{\hat{y}} = \{e \in E \hspace{1mm}|\hspace{1mm} \hat{y}_e = 1\}$ and $V_{\hat{y}} = \{i \in V \hspace{2mm}|\hspace{2mm} \exists \langle k,l \rangle \in E_{\hat{y}} : i = k \lor i = l\}$. Depending on the predictions obtained by the model, the graph may contain multiple sequences of two nodes or more. A sequence can be defined as $p = (v^p_1, v^p_2, \dots, v^p_{|p|}), v_i \in V_{\hat{y}}$, where $|p|$ is its length. Note that since the graph is undirected, the sequence can start at either of its ends. Let $E(p) = \{\{i,i+1\}\hspace{2mm}|\hspace{2mm} i \in \{1,2,\dots,|p|-1\}\}$ be the set of edges of sequence $p$ and $S(G_{\hat{y}})$ the set of all sequences in $G_{\hat{y}}$. We denote by $d_i$ the demand of node $i$ ($d_s = d_t = 0)$. The steps followed to identify and to deal with the infeasible sequences are described in Algorithm \ref{alg:check_infeasibility}.

\begin{algorithm}
\caption{Infeasibility check.}
\label{alg:check_infeasibility}
\begin{algorithmic}[1] 
 \Procedure{resolve\_infeasibility}{$G_{\hat{y}},\boldsymbol{\hat{p}}$} 
 \label{alg1_1}
    \Do    
    \label{alg1_2}
        \State \texttt{infeasibilityDetected = False} \label{alg1_3}
        \For{\texttt{each sequence $p \in S(G_{\hat{y}})$}} 
        \label{alg1_4}
        \If{$\sum_{i=1}^{|p|} d_{v^p_i}>Q$} 
        \label{alg1_5}
            \State $e = \argmin_{e \in E(p)} \hat{p}_e$
            \label{alg1_6}
            \Comment{Identify the edge with lowest probability}
            \State $E_{\hat{y}} = E_{\hat{y}} \setminus \{e\}$ 
            \label{alg1_7}
            \State $\hat{y}_e = 0$ 
            \label{alg1_8}
            \State \texttt{infeasibilityDetected = True}
            \label{alg1_9}
        \EndIf
        \label{alg1_10}
        \EndFor 
        \label{alg1_11}
    \doWhile{\texttt{infeasibilityDetected = True}}
    \label{alg1_12}
\EndProcedure

 \end{algorithmic}
\end{algorithm}

The algorithm takes as input the graph $G_{\hat{y}}$ and starts by initializing the variable \texttt{infeasibilityDetected} to \texttt{False} (Step \ref{alg1_3}). Then, it proceeds by looping over the sequences of the set $S(G_{\hat{y}})$ in Step \ref{alg1_4}. For each sequence $p$, if the total demand of the clients exceeds the capacity $Q$ of a vehicle (Step \ref{alg1_5}), the edge with the lowest probability estimate is identified (Step \ref{alg1_6}), then unfixed and removed from the set of edges $E_{\hat{y}}$ (Steps \ref{alg1_7}-\ref{alg1_8}). This procedure is repeated until there are no more infeasible sequences.

\subsubsection{\label{sec:acceleration_strat}Network reduction}

Once there are no more infeasible sequences, it is possible to make an improvement that can further accelerate the optimization, which consists in reducing the number of nodes and edges in the network. The simple (and well-known) idea is to shrink the sequences of three nodes (i.e., two edges) or more into a single edge, while making sure that the cost and the demands are updated accordingly. \\
For a sequence $p = (v^p_1, \dots, v^p_{|p|}) \in S(G_{\hat{y}})$ with $|p| \geq 3$, the goal is to remove all intermediate nodes $v^p_2, \dots, v^p_{|p|-1}$ and link directly the two ends of the sequence by adding the edge $e_{new} = \langle v^p_1,v^p_{|p|}\rangle$. Its cost is updated to the cost of the whole sequence, i.e., $c_{e_{new}} = \sum_{e\in E(p)} c_e$, and the demand of the removed nodes $\sum_{i=2}^{|p|-1} d_{v^p_i}$ is added to the demand of either $v^p_1$ or $v^p_{|p|}$. Lastly, the edge is fixed by setting the corresponding flow variable $x_{e_{new}} = 1$. By doing so, traversing the edge $e_{new}$ becomes equivalent to traversing the sequence $p$. The reduction of the network size depends on the number of sequences and their lengths, in some cases a significant reduction can be obtained, as illustrated in Figure (\ref{subfig:aggregated_edges}).

\subsubsection{\label{sec:method_summary}Method summary}

By putting all the pieces together, Algorithm \ref{alg:edge_fixing} summarizes the different steps of our method. Given the initial data $(\mathcal{P}_o, \mathcal{P}_{m}, \mathcal{S}_o)$, the algorithm starts by extracting the edge features yielding the features matrix $\mathbf{X}$ (Step \ref{alg2_1}), then obtaining the predictions $\boldsymbol{\hat{y}}$ of the ML model and the probability estimates $\boldsymbol{\hat{p}}$ in Steps \ref{alg2_2}-\ref{alg2_3}. The graph $G_{\hat{y}}$ representing the edges to fix is built (Steps \ref{alg2_4}-\ref{alg2_5}) and the procedure \texttt{resolve\_infeasibility} described in Algorithm \ref{alg:check_infeasibility} is called to check for any infeasible sequences. The algorithm proceeds by aggregating the sequences in $S(G_{\hat{y}})$ by connecting each sequence ends (Step \ref{alg2_9}), updating the cost and the demand (we added the demands of the intermediate nodes to one of the two ends) in Steps \ref{alg2_10}-\ref{alg2_11}, then removing the intermediate nodes and their adjacent edges (Steps \ref{alg2_12}-\ref{alg2_13}). Any removed node from $G_{\hat{y}}$ must be removed from the original graph also. This is achieved by updating $G$ to the induced subgraph $G[V_{\hat{y}}]$. At this point, the graph $G$ corresponds to the original complete graph minus the removed nodes from the network reduction step. Finally, the master problem is initialized and solved after adding the flow constraints for each fixed edge (Steps \ref{alg2_17}-\ref{alg2_21}).

\begin{algorithm}
\caption{ML-based edge fixing heuristic.}
\label{alg:edge_fixing}
\hspace*{\algorithmicindent} \textbf{Data:}\\ \hspace*{\algorithmicindent} $\mathcal{P}_{o}$ : Original instance\\
\hspace*{\algorithmicindent} $\mathcal{P}_{m}$ : Modified instance\\
\hspace*{\algorithmicindent} $\mathcal{S}_{o}$ : Solution of the original instance\\ 
\hspace*{\algorithmicindent} $G = (V,E)$ : Original graph
\begin{algorithmic}[1] 
 
    \State $\boldsymbol{X} \longleftarrow$ \texttt{extractEdgeFeatures($\mathcal{P}_{o}$, $\mathcal{P}_{m}$, $\mathcal{S}_{o}$)}
    \label{alg2_1}
    \State $\boldsymbol{\hat{y}} \longleftarrow$ \texttt{predict($\boldsymbol{X}$)}
    \label{alg2_2}
    \State $\boldsymbol{\hat{p}} \longleftarrow$ \texttt{predict\_proba($\boldsymbol{X}$)}
    \label{alg2_3}
    \State $V_{\hat{y}} = V, \hspace{2mm} E_{\hat{y}} = \{e \in E \hspace{1mm}|\hspace{1mm} \hat{y}_e = 1\}$
    \label{alg2_4}
    \State $G_{\hat{y}} = (V_{\hat{y}}, E_{\hat{y}})$
    \label{alg2_5}
    \State \texttt{resolve\_infeasibility}($G_{\hat{y}},\boldsymbol{\hat{p}}$)
    \label{alg2_6}
    \For{\texttt{each sequence $p = (v^p_1, \dots, v^p_{|p|}) \in S(G_{\hat{y}})$}} \hspace*{10mm} \tikzmark{right} \tikzmark{top} \label{alg2_7}
        \If{$|p| > 2$} \label{alg2_8}
            \State $e = \langle v^p_1, v^p_{|p|} \rangle$ \label{alg2_9}
            \State $c_{e} = \sum_{u \in E(p)} c_u$ \label{alg2_10}
            \State $d_{v^p_1} = d_{v^p_1} + \sum_{i=2}^{|p|-1} d_{v^p_i}$ \label{alg2_11}
            \State $V_{\hat{y}} = V_{\hat{y}} \setminus \{v^p_i \hspace{2mm}|\hspace{2mm} i \in \{2, \dots, |p|-1\}\}$ \label{alg2_12}
            \State $E_{\hat{y}} = E_{\hat{y}} \setminus \{ \langle i,j \rangle \in E \hspace{2mm}|\hspace{2mm} i \notin V_{\hat{y}} \lor j \notin V_{\hat{y}}\}$ \label{alg2_13}
        \EndIf \label{alg2_14}
    \EndFor \label{alg2_15}
    \State $G = G[V_{\hat{y}}]$ \tikzmark{bottom} \label{alg2_16}
    \State \texttt{MP $\longleftarrow$ initialize\_MP($\mathcal{P}_{m}$)} \label{alg2_17}
    \For{\texttt{each edge $e \in E_{\hat{y}}$}} \label{alg2_18}
        \State \texttt{add\_constraint(MP,"$x_e = 1$")} \label{alg2_19}
    \EndFor \label{alg2_20}
    \State \texttt{solve($\mathcal{P}_{m}, G$)} \label{alg2_21}
 \end{algorithmic}
 \AddNote{top}{bottom}{right}{Network reduction}
\end{algorithm}
\section{\label{sec:comp_experiments}Computational experiments}

This section start by describing the CVRP instances we used and the data generation process. Next, we present the details about the ML phase. Finally, the results of our heuristic method on different benchmark instances are reported. All the experiments were conducted on a Linux machine with an Xeon(R) Gold 6142 CPU @ 2.60GHz and 512GB of RAM.

\subsection{\label{sec:cvrp_instances}CVRP instances}

The instances used are based on the X benchmark dataset introduced by \citet{uchoa-x-instances}. Each instance is characterized by the following attributes:
\begin{itemize}
    \item \textbf{Depot position:} The possible values are \textbf{Central (C)} (i.e., the depot is positioned at the center of the grid), \textbf{Eccentric (E)} (i.e., the depot is positioned at the South-West corner $(0,0)$ of the grid) and \textbf{Random (R)}.
    \item \textbf{Client positioning:} The three possibilities are \textbf{Random (R)} (i.e., the clients are randomly dispersed on the grid), \textbf{Clustered (C)} (i.e., the clients are grouped in clusters) and a combination of both, referred to as \textbf{Random-Clustered (RC)}.
    \item \textbf{Demand distribution:} For the client demands, there are seven options with different intervals, and all demands are drawn uniformally from each distribution. The demands can range from: \textbf{(a)} [1-10], \textbf{(b)} [5-10], \textbf{(c)} [1-100], \textbf{(d)} [50-100], \textbf{(e)} clients located in an even quadrant have demands in the range [1-50] and [50,100] for the others, \textbf{(f)} $70\%$ to $95\%$ of the clients have a demand in the range [1,10] and [50,100] for the remaining ones and finally, \textbf{(g)} unit demands, where all the demands are equal to $1$.
    \item \textbf{Average route size:} This represents the average number of clients that can be visited by the same route, which is directly controlled by the vehicles capacity. This is computed as $n/K_{min}$ where $n$ is the number of clients and $K_{min}$ is an estimate of the minimum number of vehicles required to service all clients.
\end{itemize}
Since we are interested in modifying the demands of the clients, we exclude the instances with unit demands. Table \ref{tab:cvrp-instances} describes the instances we picked for our experiments. The instance names are written in the form ``X-n[$n_{nodes}$]-k[$n_{vehicles}$]" where $n_{nodes}$ is the number of nodes including the depot and $n_{vehicles}$ is an estimate of the number of vehicles required to service all clients. The remaining columns report the characteristics described above. 

\begin{table}[t]
  \begin{center}
    \resizebox{0.7\textwidth}{!}{
    \begin{tabular}{lcccr} 
      \toprule
      \multirow{2}{*}{\textbf{Instance name}} & \textbf{Depot} & \textbf{Client} & \textbf{Demand} & \textbf{Avg. route} \\
      & \textbf{position} & \textbf{positioning} & \textbf{distribution} & \textbf{size}\\
      \midrule
      \textbf{X-n101-k25} & R & RC & $[1-100]$ &  4.0 \\
      \textbf{X-n106-k14} & E & C & $[50-100]$ &  $7.5$ \\
      \textbf{X-n110-k13} & C & R & $[5-10]$ & 8.4 \\
      \textbf{X-n125-k30} & R & C & $Quadrant$ & 4.1 \\
      \textbf{X-n129-k18} & E & RC & $[1-10]$ & 7.1 \\
      \textbf{X-n134-k13} & R & C & $Quadrant$ & 10.2 \\
      \textbf{X-n139-k10} & C & R & $[5-10]$ & 13.8 \\
      \textbf{X-n143-k07} & E & R & $[1-100]$ & 20.3 \\
      \bottomrule
    \end{tabular}}
  \end{center}
  \caption{CVRP instances from the X benchmark instances.}
  \label{tab:cvrp-instances}
\end{table}

\subsection{\label{sec:data_generation}Data generation}

For each instance described in Table \ref{tab:cvrp-instances}, we generate a set of modified instances by randomly changing the demands of $N_{c}\%$ of the clients. The new demand of each of these clients is chosen randomly in the interval $[d_i - \Delta_d, d_i + \Delta_d],$ where $d_i$ is the original demand of the client $i$ and $\Delta_d$ a parameter controlling the interval width. In order to analyze the impact of the demand changes on the performance of the heuristic and the predictions, we used different values of $N_{c} \in \{10, 20, 30\}$. As for $\Delta_d$, since each instance has a different demand distribution, we created three classes of intervals: Small (S), Medium (M) and Large (L), controlled by the value of $\Delta_d$ as shown in Table \ref{tab:cvrp-intervals}. 

\begin{table}[t]
  \begin{center}
    \resizebox{0.4\textwidth}{!}{
    \begin{tabular}{ccccc} 
      \toprule
       \textbf{Demand} & \multicolumn{3}{c}{\textbf{$\Delta_d$ (interval size)}} \\
      \cmidrule{2-4}
       \textbf{distribution} & \hspace*{1.5mm} \textbf{S} \hspace*{1.5mm} & \hspace*{1.5mm} \textbf{M} \hspace*{1.5mm} & \hspace*{1.5mm} \textbf{L} \hspace*{1.5mm}\\
      \midrule
      $[1-100]$ & 5 & 10 & 15 \\
      $[50-100]$ & 5 & 10 & 15 \\
      $[5-10]$ & 1 & 2 & 3 \\
      $Quadrant$ & 5 & 10 & 15 \\
      $[1-10]$ & 2 & 3 & 4 \\
      \bottomrule
    \end{tabular}}
  \end{center}
  \caption{$\Delta_d$ values used depending on the demand distribution.}
  \label{tab:cvrp-intervals}
\end{table}

By combining the three different values of $N_c$ and the three interval sizes, this results in nine different scenarios $\Phi =$ \{10S, 10M, 10L, 20S, 20M, 20L, 30S, 30M, 30L\}. We then proceeded as follows. For each instance in Table \ref{tab:cvrp-instances} and each scenario $\phi \in \Phi$, $100$ modified instances are generated, where $95$ of them are used for the ML phase (i.e., training, parameters tuning, etc.) and the remaining $5$ for the optimization phase (i.e., when the ML model is incorporated in the CG algorithm). We therefore consider the eight instances of Table \ref{tab:cvrp-instances} as the original instances (i.e., $\mathcal{P}_o^1, \dots, \mathcal{P}_o^8$) and, for each scenario $\phi$, we generate $100$ modified versions (i.e., $(\mathcal{P}_{m}^{i})_{\phi}^{1}, \dots, (\mathcal{P}_{m}^{i})_{\phi}^{100}, i \in \{1,2,\dots,8\}, \phi \in \Phi$). One ML-model is trained for each instance and for each scenario. The idea of training a model for each instance comes from the assumption that we have a specific instance that we solve repeatedly. Therefore we want a specific model for that particular instance. This makes a total of $8$ instances $\times \hspace{1mm}9$ scenarios  $\times \hspace{1mm}100 = 7,200$ modified instances generated (and $8 \times 9 = 72$ ML models). Given the large number of instances, instead of using an exact B\&P method for solving and collecting solutions, we opted to use the FILO heuristic of \citet{filo}(see Section \ref{sec:cvrp-heuristic}). For each of the $7,200$ instances, we run the heuristic using $10$ different random seed values for 1,000,000 iterations and the solution with the smallest cost is retained.

Once the solutions of the original and modified instances are obtained, for each instance $i \in \{1,\dots,8\}$ and scenario $\phi \in \Phi$, the solutions are grouped in a set of tuples, i.e., $\{(\mathcal{P}_o^i, \mathcal{S}_o^i, (\mathcal{P}_m^i)_{\phi}^1, (\mathcal{S}_m^i)_{\phi}^1),\allowbreak \dots,\allowbreak (\mathcal{P}_o^i, \mathcal{S}_o^i, (\mathcal{P}_m^i)_{\phi}^{100}, (\mathcal{S}_m^i)_{\phi}^{100})\}$ that are used to extract the edge features and labels as detailed in Section \ref{sec:data-collect}.

\subsection{\label{sec:ml_phase}Machine learning phase}

Before starting the training, common practices in ML are followed, starting by a pre-processing phase that consists of scaling and normalizing the data. The dataset (i.e., the data from the 95 instances) is then split into a training set, a validation set and a test set. The goal of the validation set is to tune the different hyperparameters, whereas the purpose of the test set is to compare different classification algorithms, such as logistic regression, K-nearest neighbors, random forest, artificial neural network (ANN), etc. According to the results obtained on the test set, the ANN model is the overall most robust model with accuracies ranging from $70\%$ to $88\%$ depending on the instance and the scenario $\phi \in \Phi$ (more detailed results about the models performance are reported in the next section). The hyperparameter values used during the training of the ANN models are described in Table \ref{tab:ann-hyperpameters}.

\begin{table}[t]
  \begin{center}
    
    \begin{tabular}{lc} 
      \toprule
      \textbf{Hyperparameter} & \textbf{Value} \\
      \midrule
       Learning rate & $10^{-3}$\\
       Number of epochs & 1000\\
       Epoch size & 64\\
       Batch size & 32\\
       NN Architecture & $32\times32\times32\times1$\\
       Activation function & ReLU\\
       Output function & Sigmoid\\
       Optimizer & Adam\\
       Class weights & Balanced\\
      \bottomrule
    \end{tabular}
  \end{center}
  \caption{Hyperparameters values of the ANN models.}
  \label{tab:ann-hyperpameters}
\end{table}


\subsection{\label{sec:heuristic_results}Optimization phase}

In this section, we present the results obtained by our edge-fixing heuristic. Since a ML model is trained for each original instance and scenario, we report the ML model performance in this section as well. Tables \ref{tab:cvrp_summary_1} to \ref{tab:cvrp_summary_3} summarize the results obtained on the test instances. Each row corresponds to the average values obtained on the 5 test instances (for each original instance), whereas the details of each individual instance can be found in Appendix \ref{sec:appendix_A}. There is one table for each $N_c$ value (i.e., 10, 20 and 30 for Tables \ref{tab:cvrp_summary_1}, \ref{tab:cvrp_summary_2} and \ref{tab:cvrp_summary_3}, respectively). The first column corresponds to the interval used when changing the demands (i.e., S, M and L, respectively), followed by the name of the original instance $\mathcal{P}_o$. Next, the average cost of the best solutions $\mathcal{S}_m$ obtained by the FILO heuristic (over the 10 executions of the heuristic with different seeds). In the fourth column, we report the average similarity between the solution of the original instance and the solution of the modified instances computed by the following formula:

\begin{align}
sim(\mathcal{S}_{o},\mathcal{S}_{m}) = \frac{|E(\mathcal{S}_{o}) \cap E(\mathcal{S}_{m})|}{|E(\mathcal{S}_{o})|}.
\end{align}

Notice that this similarity also matches the percentage of edges to fix (i.e., with the label $1$). The next three columns summarize the performance of the ML model and show the following metrics: the True Negative Rate (\textbf{TNR}) corresponds to the percentage of edges that should not be fixed and that are predicted accurately; the True Positive Rate (\textbf{TPR}) is equal to the percentage of edges that should be fixed and that are predicted correctly; and the (balanced) accuracy is the mean of the two previous columns. By fixing the edges suggested by the ML model and solving the instance using VRPSolver, we obtain the results shown under the heading ``Edge-fixing": the average costs of the solutions as well as the computing times in seconds. To evaluate the quality of this solution, we report the average gap that compares the cost of the solution of our method with that of the solution $\mathcal{S}_m$ of the FILO heuristic (third column). The next three columns provide the average cost of the solutions computed by the exact algorithm in VRPSolver, the average computing time and the ratio with respect to the computing time of our method (i.e., the computing time of the exact B\&P algorithm divided by the computing time of our edge-fixing method). For the exact algorithm, a time limit of five hours is set. Empty values mean that the exact B\&P failed to find an optimal solution in the time limit for one or more of the 5 instances, refer to Appendix \ref{sec:appendix_A} for additional details. 

A natural question with respect to the proposed approach is if there is value in trying to learn the difference between solutions of slightly modified instances. Indeed, ML-based CVRP heuristics should be in a very favorable position in our computational setting because the instances are from a specific data distribution and of the same size, i.e., no generalization is needed. In order to give a tentative answer to this question, we compare the performance of our method with another ML method for the CVRP not designed for reoptimization purposes but with excellent overall performance on recent benchmarks. We chose the \textit{Dual-Aspect Collaborative Transformer} algorithm (DACT - \citet{dact-ml}) that is considered one of the most effective methods for solving the CVRP to-date in terms of solution quality according to the recent survey \citet{cvrp-ml-survey}. Unlike our approach, DACT is a reinforcement learning method that learns an improvement heuristic, which means that it starts with an initial solution and tries to improve it in an iterative way. In this method, the learner (i.e., the agent) learns to identify a pair of nodes on which to apply a pairwise operator (e.g., 2-opt, swap, insert) and is rewarded when a better solution is found. The authors report good results that outperform several other learning methods on the X benchmark instances (same instances we are using) in a reasonable computing time. In our case, we proceeded by using the same pretrained model that the authors used in their paper but with additional training on the 8 instances we focus on (see Table \ref{tab:cvrp-instances}), in addition to using the same parameters of their best performing model. The results obtained are reported in the last two columns of Tables \ref{tab:cvrp_summary_1} to \ref{tab:cvrp_summary_3}, representing respectively the average cost of the solutions and the gap with respect to the FILO heuristic (just like Edge-Fixing).

\begin{table}[t]
  \begin{center}
    \resizebox{\textwidth}{!}{
    \begin{tabular}{ccrrrrrrrrrrrrr} 
      \toprule
      \multirow{2}{*}{\textbf{Interval}} &
      \multirow{2}{*}{$\mathcal{P}_{o}$} & \multirow{2}{*}{\textbf{$\mathcal{S}_m$ cost}} &  \multirow{2}{*}{\textbf{$sim(\mathcal{S}_{o},\mathcal{S}_{m})$}} & \multicolumn{3}{c}{\textbf{ML model metrics}} & \multicolumn{3}{c}{\textbf{Edge-fixing}} & \multicolumn{3}{c}{\textbf{Exact B\&P}} & \multicolumn{2}{c}{\textbf{DACT}}\\
      \cmidrule(rl){5-7} \cmidrule(rl){8-10} \cmidrule(rl){11-13} \cmidrule(rl){14-15}
       & & & & \textbf{TNR} & \textbf{TPR} & \textbf{Accuracy} & \textbf{Cost} & \textbf{Time (s)} & \textbf{Gap} & \textbf{Cost} & \textbf{Time (s)} & \textbf{Ratio} & \textbf{Cost} & \textbf{Gap}\\
      \midrule

      \multirow{8}{*}{\textbf{Small}} & X-n101-k25 & 27637 & 84\% & 71\% & 70\% & 70\% & 27718 & 13 & 0.29\% & 27635 & 211  & 22.4 & 28184 & 1.98\%\\
      & X-n106-k14 & 26376 & 85\% & 89\% & 60\% & 75\% & 26436 & 25 & 0.23\% & 26376 & 962 & 61.3 & 26897 & 1.98\%\\
      & X-n110-k13 & 14987 & 93\% & 100\% & 66\% & 83\% & 14987 & 8 & 0.00\% & 14987 & 265 & 36.7 & 15159 & 1.15\% \\
        & X-n125-k30 & 55613 & 63\% & 76\% & 73\% & 75\% & 55683 & 235 & 0.13\% & - & - & - & 58581 & 5.34\% \\
        & X-n129-k18 & 28765 & 62\% & 78\% & 72\% & 75\% & 28982 & 78 & 0.75\% & 28765 & 3479 & 53.4 & 29697 & 3.24\% \\
        & X-n134-k13 & 10888 & 80\% & 89\% & 54\% & 72\% & 10917 & 77 & 0.27\% & - & - & - & 11284 & 3.64\%\\
        & X-n139-k10 & 13590 & 85\% & 92\% & 81\% & 86\% & 13599 & 26 & 0.06\% & - & - & - & 13846 & 1.88\% \\
        & X-n143-k07 & 15722 & 81\% & 87\% & 88\% & 88\% & 15726 & 54 & 0.02\% & - & - & - & 16245 & 3.32\% \\
        \multicolumn{2}{l}{\textbf{Average}} & \textbf{23886} & \textbf{79\%} & \textbf{85\%} & \textbf{71\%} & \textbf{78\%} & \textbf{24256} & \textbf{65} & \textbf{0.22\%} & \textbf{-} & \textbf{-} & \textbf{-} & \textbf{24987} & \textbf{2.82\%}\\
        \midrule
        \multirow{8}{*}{\textbf{Medium}} & X-n101-k25 & 27606 & 83\% & 77\% & 63\% & 70\% & 27736 & 18 & 0.47\% & 27606 & 324 & 20.4 & 28298 & 2.51\% \\
        & X-n106-k14 & 26358 & 66\% & 95\% & 73\% & 84\% & 26380 & 14 & 0.08\% & 26358 & 355 & 23.9 & 26871 & 1.95\% \\
        & X-n110-k13 & 14971 & 87\% & 98\% & 67\% & 82\% & 14993 & 12 & 0.15\% & 14969 & 283 & 24.5 & 15137 & 1.11\% \\
        & X-n125-k30 & 55713 & 60\% & 74\% & 73\% & 74\% & 55758 & 198 & 0.08\% & 55655 & 5156 & 58.4 & 58735 & 5.42\% \\
        & X-n129-k18 & 28862 & 60\% & 74\% & 75\% & 75\% & 29124 & 105 & 0.91\% & - & - & - & 29801 & 3.26\% \\
        & X-n134-k13 & 10888 & 63\% & 74\% & 73\% & 74\% & 11024 & 320 & 1.25\% & - & - & - & 11304 & 3.82\% \\
        & X-n139-k10 & 13601 & 90\% & 85\% & 76\% & 81\% & 13608 & 27 & 0.05\% & - & - & - & 13863 & 1.92\% \\
        & X-n143-k07 & 15707 & 85\% & 97\% & 79\% & 88\% & 15710 & 63 & 0.02\% & - & - & - & 16200 & 3.14\% \\
        \multicolumn{2}{l}{\textbf{Average}} & \textbf{24213} & \textbf{74\%} & \textbf{84\%} & \textbf{72\%} & \textbf{78\%} & \textbf{24292} & \textbf{95} & \textbf{0.38\%} & \textbf{-} & \textbf{-} & \textbf{-} & \textbf{25026} & \textbf{2.89\%}\\
        \midrule
        \multirow{8}{*}{\textbf{Large}} & X-n101-k25 & 27651 & 70\% & 63\% & 77\% & 70\% & 28125 & 122 & 1.71\% & 27648 & 399 & 7.9 & 28142 & 1.98\% \\
        & X-n106-k14  & 26412 & 65\% & 84\% & 74\% & 79\% & 26557 & 250 & 0.54\% & - & - & - & 26846 & 1.64\% \\
        & X-n110-k13  & 15030 & 75\% & 86\% & 73\% & 80\% & 15107 & 35 & 0.51\% & 15030 & 2283 & 60.7 & 15187 & 1.04\% \\
        & X-n125-k30  & 55733 & 55\% & 81\% & 75\% & 78\% & 55825 & 422 & 0.16\% & - & - & - & 58500 & 4.97\% \\
        & X-n129-k18  & 28755 & 62\% & 78\% & 73\% & 75\% & 28972 & 171 & 0.76\% & 28748 & 2457 & 22.5 & 29546 & 2.75\% \\
        & X-n134-k13 & 10908 & 69\% & 72\% & 68\% & 70\% & 10908 & 164 & 0.66\% & - & - & - & 11267 & 3.29\% \\
        & X-n139-k10  & 13600 & 83\% & 85\% & 75\% & 80\% & 13635 & 47 & 0.26\% & - & - & - & 13850 & 1.84\% \\
        & X-n143-k07  & 15716 & 88\% & 97\% & 73\% & 85\% & 15717 & 62 & 0.00\% & - & - & - & 16265 & 3.49\%\\
        \multicolumn{2}{l}{\textbf{Average}} & \textbf{24226} & \textbf{71\%} & \textbf{81\%} & \textbf{73\%} & \textbf{77\%} & \textbf{24365} & \textbf{159} & \textbf{0.58\%} & \textbf{-} & \textbf{-} & \textbf{-} & \textbf{24950} & \textbf{2.62\%}\\
        \bottomrule
    \end{tabular}}
  \end{center}
  \caption{Average results for scenarios with $N_c = 10$.}
  \label{tab:cvrp_summary_1}
\end{table}

\begin{table}[t]
  \begin{center}
    \resizebox{\textwidth}{!}{
    \begin{tabular}{ccrrrrrrrrrrrrr} 
      \toprule
      \multirow{2}{*}{\textbf{Interval}} &
      \multirow{2}{*}{$\mathcal{P}_{o}$} & \multirow{2}{*}{\textbf{$\mathcal{S}_m$ cost}} &  \multirow{2}{*}{\textbf{$sim(\mathcal{S}_{o},\mathcal{S}_{m})$}} & \multicolumn{3}{c}{\textbf{ML model metrics}} & \multicolumn{3}{c}{\textbf{Edge-fixing}} & \multicolumn{3}{c}{\textbf{Exact B\&P}} & \multicolumn{2}{c}{\textbf{DACT}}\\
      \cmidrule(rl){5-7} \cmidrule(rl){8-10} \cmidrule(rl){11-13} \cmidrule(rl){14-15}
       & & & & \textbf{TNR} & \textbf{TPR} & \textbf{Accuracy} & \textbf{Cost} & \textbf{Time (s)} & \textbf{Gap} & \textbf{Cost} & \textbf{Time (s)} & \textbf{Ratio} & \textbf{Cost} & \textbf{Gap}\\
      \midrule
      
        \multirow{8}{*}{\textbf{Small}} & X-n101-k25 & 27486 & 83\% & 73\% & 75\% & 74\% & 27628 & 29 & 0.51\% & 27486 & 231 & 8.6 & 28220 & 2.67\% \\
        & X-n106-k14  & 26338 & 72\% & 90\% & 71\% & 81\% & 26401 & 30 & 0.24\% & 26336 & 1147 & 88.5 & 26845 & 1.93\% \\
        & X-n110-k13  & 14980 & 77\% & 86\% & 74\% & 80\% & 15024 & 9 & 0.29\% & 14980 & 429 & 55.5 & 15152 & 1.15\% \\
        & X-n125-k30  & 55493 & 63\% & 85\% & 68\% & 77\% & 55509 & 233 & 0.03\% & - & - & - & 58372 & 5.19\% \\
        & X-n129-k18 & 29020 & 56\% & 71\% & 78\% & 74\% & 29408 & 143 & 1.33\% & 29009 & 5791 & 76.5 & 29703 & 2.35\% \\
        & X-n134-k13  & 10909 & 73\% & 85\% & 69\% & 77\% & 10943 & 150 & 0.32\% & - & - & - & 11323 & 3.80\% \\
        & X-n139-k10  & 13613 & 87\% & 91\% & 81\% & 86\% & 13616 & 103 & 0.02\% & - & - & - & 13870 & 1.89\% \\
        & X-n143-k07  & 15715 & 86\% & 90\% & 86\% & 88\% & 15748 & 30 & 0.21\% & - & - & - & 16261 & 3.47\% \\
        \multicolumn{2}{l}{\textbf{Average}} & \textbf{24194} & \textbf{75\%} & \textbf{84\%} & \textbf{75\%} & \textbf{80\%} & \textbf{24284} & \textbf{91} & \textbf{0.37\%} & \textbf{-} & \textbf{-} & \textbf{-} & \textbf{24968} & \textbf{2.81\%}\\
        \midrule
        \multirow{8}{*}{\textbf{Medium}} & X-n101-k25 & 27512 & 68\% & 74\% & 72\% & 73\% & 27694 & 32 & 0.66\% & 27512 & 518 & 19.4 & 28076 & 2.05\% \\
        & X-n106-k14 & 26306 & 64\% & 87\% & 71\% & 79\% & 26400 & 20 & 0.36\% & 26306 & 5257 & 238.2 & 26785 & 1.72\% \\
        & X-n110-k13 & 14983 & 73\% & 86\% & 77\% & 82\% & 15078 & 27 & 0.64\% & 14983 & 548 & 30.3 & 15217 & 1.56\% \\
        & X-n125-k30 & 55503 & 54\% & 82\% & 74\% & 78\% & 55579 & 233 & 0.14\% & - & - & - & 58528 & 5.45\% \\
        & X-n129-k18 & 29171 & 57\% & 73\% & 76\% & 74\% & 29655 & 315 & 1.66\% & - & - & - & 30001 & 2.85\% \\
        & X-n134-k13 & 10877 & 55\% & 82\% & 80\% & 81\% & 10956 & 127 & 0.73\% & - & - & - & 11258 & 3.50\% \\
        & X-n139-k10 & 13605 & 72\% & 87\% & 80\% & 84\% & 13631 & 205 & 0.20\% & - & - & - & 13855 & 1.84\% \\
        & X-n143-k07 & 15708 & 81\% & 80\% & 83\% & 82\% & 15745 & 49 & 0.24\% & - & - & - & 16233 & 3.35\% \\
        \multicolumn{2}{l}{\textbf{Average}} & \textbf{24208} & \textbf{65\%} & \textbf{81\%} & \textbf{77\%} & \textbf{79\%} & \textbf{24342} & \textbf{126} & \textbf{0.58\%} & \textbf{-} & \textbf{-} & \textbf{-} & \textbf{24991} & \textbf{2.79\%}\\
        \midrule
        \multirow{8}{*}{\textbf{Large}} & X-n101-k25 & 27645 & 59\% & 70\% & 73\% & 72\% & 27871 & 56 & 0.81\% & 27645 & 466 & 30.8 & 28244 & 2.17\% \\
        & X-n106-k14 & 26413 & 61\% & 79\% & 74\% & 77\% & 26555 & 100 & 0.54\% & - & - & - & 26834 & 1.59\% \\
        & X-n110-k13 & 15034 & 68\% & 81\% & 79\% & 80\% & 15161 & 23 & 0.84\% & 15034 & 487 & 25.4 & 15216 & 1.21\% \\
        & X-n125-k30 & 55958 & 58\% & 79\% & 70\% & 74\% & 56168 & 433 & 0.37\% & - & - & - & 58629 & 4.77\% \\
        & X-n129-k18 & 29123 & 55\% & 77\% & 74\% & 75\% & 29462 & 386 & 1.16\% & 29092 & 3106 & 16.0 & 29848 & 2.49\% \\
        & X-n134-k13 & 10909 & 55\% & 82\% & 80\% & 81\% & 11051 & 254 & 1.31\% & - & - & - & 11317 & 3.75\% \\
        & X-n139-k10 & 13585 & 73\% & 88\% & 77\% & 83\% & 13619 & 139 & 0.25\% & - & - & - & 13788 & 1.50\% \\
        & X-n143-k07 & 15743 & 84\% & 79\% & 79\% & 79\% & 15775 & 41 & 0.20\% & - & - & - & 16332 & 3.74\%\\
        \multicolumn{2}{l}{\textbf{Average}} & \textbf{24301} & \textbf{64\%} & \textbf{79\%} & \textbf{76\%} & \textbf{78\%} & \textbf{24458} & \textbf{179} & \textbf{0.69\%} & \textbf{-} & \textbf{-} & \textbf{-} & \textbf{25026} & \textbf{2.65\%}\\
        \midrule
      
    \end{tabular}}
  \end{center}
  \caption{Average results for scenarios with $N_c = 20$.}
  \label{tab:cvrp_summary_2}
\end{table}

\begin{table}[t]
  \begin{center}
    \resizebox{\textwidth}{!}{
    \begin{tabular}{ccrrrrrrrrrrrrr} 
      \toprule
      \multirow{2}{*}{\textbf{Interval}} &
      \multirow{2}{*}{$\mathcal{P}_{o}$} & \multirow{2}{*}{\textbf{$\mathcal{S}_m$ cost}} &  \multirow{2}{*}{\textbf{$sim(\mathcal{S}_{o},\mathcal{S}_{m})$}} & \multicolumn{3}{c}{\textbf{ML model metrics}} & \multicolumn{3}{c}{\textbf{Edge-fixing}} & \multicolumn{3}{c}{\textbf{Exact B\&P}} & \multicolumn{2}{c}{\textbf{DACT}}\\
      \cmidrule(rl){5-7} \cmidrule(rl){8-10} \cmidrule(rl){11-13} \cmidrule(rl){14-15}
       & & & & \textbf{TNR} & \textbf{TPR} & \textbf{Accuracy} & \textbf{Cost} & \textbf{Time (s)} & \textbf{Gap} & \textbf{Cost} & \textbf{Time (s)} & \textbf{Ratio} & \textbf{Cost} & \textbf{Gap}\\
      \midrule
      
      \multirow{8}{*}{\textbf{Small}} & X-n101-k25 & 27562 & 75\% & 80\% & 70\% & 75\% & 27669 & 37 & 0.39\% & 27562 & 178 & 5.9 & 28219 & 2.38\% \\
        & X-n106-k14 & 26383 & 65\% & 89\% & 69\% & 79\% & 26438 & 148 & 0.21\% & 26378 & 1577 & 25.9 & 26833 & 1.70\% \\
        & X-n110-k13 & 15005 & 74\% & 92\% & 77\% & 84\% & 15083 & 32 & 0.52\% & 15005 & 555 & 21.8 & 15157 & 1.01\% \\
        & X-n125-k30 & 55776 & 53\% & 77\% & 73\% & 75\% & 55834 & 259 & 0.10\% & - & - & - & 58601 & 5.06\% \\
        & X-n129-k18 & 29414 & 58\% & 80\% & 74\% & 77\% & 29778 & 237 & 1.24\% & 29405 & 3742 & 24.8 & 30086 & 2.28\% \\
        & X-n134-k13 & 10925 & 77\% & 90\% & 62\% & 76\% & 10952 & 140 & 0.25\% & - & - & - & 11361 & 4.00\% \\
        & X-n139-k10 & 13622 & 83\% & 91\% & 73\% & 82\% & 13692 & 154 & 0.51\% & - & - & - & 13829 & 1.52\% \\
        & X-n143-k07 & 15754 & 83\% & 79\% & 87\% & 83\% & 15822 & 56 & 0.43\% & - & - & - & 16328 & 3.65\% \\
        \multicolumn{2}{l}{\textbf{Average}} & \textbf{24305} & \textbf{71\%} & \textbf{85\%} & \textbf{73\%} & \textbf{79\%} & \textbf{24409} & \textbf{133} & \textbf{0.46\%} & \textbf{-} & \textbf{-} & \textbf{-} & \textbf{25052} & \textbf{2.70\%}\\
        \midrule
        
        \multirow{8}{*}{\textbf{Medium}} & X-n101-k25 & 27654 & 63\% & 68\% & 73\% & 71\% & 27804 & 51 & 0.54\% & 27651 & 243 & 24.8 & 28277 & 2.25\% \\
        & X-n106-k14 & 26437 & 59\% & 75\% & 70\% & 73\% & 26615 & 247 & 0.67\% & - & - & - & 26780 & 1.30\% \\
        & X-n110-k13 & 15058 & 77\% & 87\% & 72\% & 79\% & 15095 & 14 & 0.24\% & 15058 & 269 & 33.6 & 15292 & 1.55\% \\
        & X-n125-k30 & 55925 & 50\% & 82\% & 78\% & 80\% & 56009 & 451 & 0.15\% & - & - & - & 58859 & 5.25\% \\
        & X-n129-k18 & 29221 & 54\% & 79\% & 75\% & 77\% & 29698 & 137 & 1.63\% & 29215 & 3574 & 34.8 & 29849 & 2.16\% \\
        & X-n134-k13 & 10926 & 56\% & 83\% & 80\% & 81\% & 11048 & 129 & 1.12\% & - & - & - & 11301 & 3.43\% \\
        & X-n139-k10 & 13640 & 75\% & 89\% & 78\% & 83\% & 13708 & 256 & 0.49\% & - & - & - & 13926 & 2.09\% \\
        & X-n143-k07 & 15782 & 80\% & 78\% & 85\% & 82\% & 15870 & 49 & 0.55\% & - & - & - & 16423 & 4.06\% \\
        \multicolumn{2}{l}{\textbf{Average}} & \textbf{24330} & \textbf{64\%} & \textbf{80\%} & \textbf{76\%} & \textbf{78\%} & \textbf{24481} & \textbf{167} & \textbf{0.67\%} & \textbf{-} & \textbf{-} & \textbf{-} & \textbf{25088} & \textbf{2.76\%}\\
        \midrule
      
        \multirow{8}{*}{\textbf{Large}} & X-n101-k25 & 27617 & 66\% & 76\% & 72\% & 74\% & 27929 & 185 & 1.12\% & 27617 & 339 & 19.4 & 28207 & 2.13\% \\
        & X-n106-k14 & 26440 & 54\% & 79\% & 68\% & 74\% & 26637 & 275 & 0.75\% & - & - & - & 26913 & 1.79\% \\
        & X-n110-k13 & 15090 & 66\% & 82\% & 72\% & 77\% & 15211 & 88 & 0.80\% & 15090 & 1120 & 18.7 & 15318 & 1.51\% \\
        & X-n125-k30 & 56228 & 52\% & 84\% & 71\% & 78\% & 56221 & 475 & -0.01\% & - & - & - & 58725 & 4.44\% \\
        & X-n129-k18 & 29674 & 58\% & 75\% & 73\% & 74\% & 30179 & 75 & 1.69\% & 29670 & 4914 & 69.9 & 30554 & 2.97\% \\
        & X-n134-k13 & 10950 & 53\% & 83\% & 82\% & 82\% & 11029 & 145 & 0.72\% & - & - & - & 11423 & 4.31\% \\
        & X-n139-k10 & 13616 & 73\% & 91\% & 76\% & 84\% & 13635 & 128 & 0.14\% & - & - & - & 13868 & 1.85\% \\
        & X-n143-k07 & 15884 & 79\% & 82\% & 83\% & 83\% & 15941 & 313 & 0.36\% & - & - & - & 16326 & 2.79\%\\
        \multicolumn{2}{l}{\textbf{Average}} & \textbf{24437} & \textbf{63\%} & \textbf{82\%} & \textbf{75\%} & \textbf{78\%} & \textbf{24598} & \textbf{211} & \textbf{0.70\%} & \textbf{-} & \textbf{-} & \textbf{-} & \textbf{25167} & \textbf{2.72\%}\\
        \midrule
    \end{tabular}}
  \end{center}
  \caption{Average results for scenarios with $N_c = 30$.}
  \label{tab:cvrp_summary_3}
\end{table}

According to the results reported in Tables \ref{tab:cvrp_summary_1} to \ref{tab:cvrp_summary_3}, we can notice that on average the similarity decreases by increasing the number of changes we make on the demands (controlled by the parameter $N_c$ and the intervals), which is expected. A high similarity of 79\% is noticed on average on the instances with 10\% change and small interval, while the instances with 30\% change and large interval have an average similarity of 63\%. At the instance level, it seems that the dispersion of the clients as well as the length of the routes can have an impact on the similarity of the solutions. For example, the instances with randomly dispersed clients and medium to long routes (e.g., X-n110-k13, X-n139-k10, X-n143-k07) tend to have a higher similarity compared to instances with clustered clients (e.g., X-n106-k14 and X-n125-k30).

Before evaluating the performance of the model, we would like to highlight the impact that a bad prediction may have in the obtained solution. If some edges are not fixed when they should be (false negatives), it can affect the computing time but the solver can still include them in the final solution. However, the edges with label 0 can have a more serious impact, since they can affect the quality of the solution if predicted inaccurately (false positives), which potentially leads to a higher gap. As previously mentioned, the similarity is also the percentage of edges to be fixed. A high similarity means that most of the edges in the solution of the original instance are labeled $1$, which also implies that the ML model tends to make fewer ``significant" errors since there are not many edges with label $0$. If we focus on the solution quality, it is possible to give a higher weight to the edges with label 0, which can lead to a high TNR (thus reducing the false positives) and probably a lower TPR, meaning that less edges will be fixed resulting in a higher computing time. Conversely, if we do the opposite and assign a higher weight to the edges with label 1, we will likely fix more edges and achieve a lower computing time but at the expense of bad-quality solutions. In our case, we use a balanced weight between the two classes as shown in Table \ref{tab:ann-hyperpameters}. Depending on the desired goal, one seeks to find a compromise between quality and computing time. 

From the ML metrics obtained we can observe an average accuracy ranging from $78\%$ to $80\%$, but if we look closer at the individual instances, we can notice some variance and a slight positive correlation between the accuracy and the similarity. By fixing the edges proposed by the model and comparing the solution with the one of the FILO heuristic, we obtain an average gap ranging from $0.22\%$ to $0.70\%$. We can also notice that the gap follows the same tendency as the similarity, i.e., more changes of the demands imply less similarity, more possible misclassifications and thus a higher gap. Therefore, the instances with high similarity (e.g., instances with dispersed clients and long routes) generally have a lower gap compared to the others. The detailed results in Appendix \ref{sec:appendix_A} indicate that the gap of the individual instances can vary between $0\%$ and $1.71\%$. A perfect solution is obtained (i.e., same solution as the FILO heuristic, thus a $0\%$ gap) when the TNR is $100\%$. In a few rare cases, the solver was able to find a better solution than the heuristic even after fixing parts of the solution, such as X-n125-k30\_10M\_4 in Table \ref{tab:cvrp_results_2} and X-n125-k30\_30L\_3 in Table \ref{tab:cvrp_results_9}. Generally, the FILO heuristic finds very good solutions (especially when considering 10 runs): it succeeds in finding the optimal solution in many cases when compared to the available exact B\&P results. We believe that the gap reported between the edge-fixing and the FILO heuristic must not be far from the one with the exact B\&P algorithm.

In terms of computing time, the edge fixing method takes only a few minutes or even a few seconds in some cases to find the optimal solution of the modified problem (that with fixed variables according to the ML prediction), whereas the exact B\&P algorithm that solves the (original) problem from scratch can take hours of computation to find an optimal solution, and in many cases no solution is found after the time limit especially for large instances. Apparently, the current version of VRPSolver does not have strong heuristics to produce feasible solutions. On the one side, this is not an issue for our computational investigation because we use VRPSolver as well (and we are definitely not using it for comparison). On the other side, this is somehow reinforcing the fact that, after our proposed fixing, the instances are much easier and even without good heuristics VRPSolver is very effective on them. Note that in our experiments, the FILO heuristic with 1,000,000 iterations takes about 8 to 13 minutes per execution depending on the instance size, not to forget that we executed the heuristic 10 different times with different seeds and retained the best solution. Overall, the computing time of the edge-fixing heuristic remains lower than a single execution of FILO and way lower than an exact B\&P algorithm.

Concerning the DACT method, we can notice an average gap of $2.7\%$. The method performs better on some instances than on others with gaps ranging from $0.03\%$ up to $5.72\%$ according to the detailed results in Appendix \ref{sec:appendix_A}. Since the method is iterative and performs a certain number of steps at the inference phase (more precisely 10,000 steps as described in \citet{dact-ml}), its computing time is not negligible and can range from 10 to 15 minutes on average depending on the instance size. Therefore, the edge-fixing method outperforms DACT both in terms of computing time and solution quality.

\section{\label{sec:conclusion}Conclusion}

In this paper, we proposed a ML-based heuristic for the reoptimization of repeatedly solved problems after minor changes in the problem data. More precisely, we put ourselves in the context of a company that solves CVRPs on a daily basis where the locations of the clients are the same but slightly different demands occur. Given that there can be a great similarity between the solutions, the goal was to exploit the ones obtained from previous executions (not necessarily optimal ones) in order to speed up the reoptimization of future instances. The aim was to predict and fix the edges that have a high chance of remaining in the solution after a change of the demands.

Following a supervised learning approach, we trained neural network models on the data collected from a recent heuristic for the CVRP called FILO. The models achieved an average accuracy of 78\% on the test instances. By incorporating the predictions in our edge-fixing method, we were able to find solutions in reasonable computing times with gaps ranging from $0\%$ to $1.71\%$ (with an average of $0.51\%$) when compared to the FILO heuristic we learned from. These gaps are also lower compared to other recent ML methods proposed in the literature such as the DACT method used to compute the entire CVRP solution. We also obtained an acceleration effect, considerably reducing the size of the network and allowing an even faster reoptimization. 

Future work can seek some improvements, for example by the exploration of more complex learning models that can take advantage of the graph structure of the problem or even sequence models since we are dealing with routes. This may potentially lead to a better accuracy and therefore better gaps and lower computing times.

\newpage
\clearpage
\bibliography{main} 

\newpage
\appendix

\section{\label{sec:appendix_A}Detailed results of the experiments}

This section provides detailed computational results that are complementary to those presented in Section \ref{sec:comp_experiments}. The results are grouped by scenario ($N_c$ and interval) and reported in Tables \ref{tab:cvrp_results_1}-\ref{tab:cvrp_results_9}. The same information as in Tables \ref{tab:cvrp_summary_1} to \ref{tab:cvrp_summary_3} is reported, with the only difference that each row represents an individual instance. Also, an additional column has been added to indicate the name suffix of the modified instance $\mathcal{P}_{m}$ (second column).

\begin{table}[H]
  \begin{center}
    \resizebox{\textwidth}{!}{
\begin{tabular}{ccrrrrrrrrrrrrr} 
      \toprule
      \multirow{2}{*}{$\mathcal{P}_{o}$} &
      \multirow{2}{*}{$\mathcal{P}_{m}$} & \multirow{2}{*}{\textbf{$\mathcal{S}_m$ cost}} &  \multirow{2}{*}{\textbf{$sim(\mathcal{S}_{o},\mathcal{S}_{m})$}} & \multicolumn{3}{c}{\textbf{ML model metrics}} & \multicolumn{3}{c}{\textbf{Edge-fixing}} & \multicolumn{3}{c}{\textbf{Exact BaP}} & \multicolumn{2}{c}{\textbf{DACT}}\\
      \cmidrule(rl){5-7} \cmidrule(rl){8-10} \cmidrule(rl){11-13} \cmidrule(rl){14-15}
       & & & & \textbf{TNR} & \textbf{TPR} & \textbf{Accuracy} & \textbf{Cost} & \textbf{Time (s)} & \textbf{Gap} & \textbf{Cost} & \textbf{Time (s)} & \textbf{Ratio} & \textbf{Cost} & \textbf{Gap}\\
      \midrule
      
      \multirow{5}{*}{\textbf{X-n101-k25}} & \textbf{\_10S\_1} & 27538  & 84\% & 55\% & 67\% & 61\% & 27733 & 26 & 0.71\% & 27538 & 226 & 8.7 & 27943 & 1.47\%\\
      & \textbf{\_10S\_2} & 27775  & 71\% & 61\% & 73\% & 67\% & 27849 & 14 & 0.27\% & 27764 & 199 & 14.2 & 28306 & 1.91\%\\ 
      & \textbf{\_10S\_3} & 27591  & 100\% & 100\% & 63\% & 82\% & 27591 & 7 & 0.00\% & 27591 & 180 & 25.7 & 28287 & 2.52\% \\ 
      & \textbf{\_10S\_4} & 27703  & 73\% & 60\% & 76\% & 68\% & 27811 & 14 & 0.39\% & 27703 & 276 & 19.7 & 27966 & 0.95\%\\ 
      & \textbf{\_10S\_5} & 27580  & 91\% & 77\% & 70\% & 74\% & 27605 & 4 & 0.09\% & 27580 & 175 & 43.8 & 28417 & 3.03\%\\
      \textbf{Average}  & - & \textbf{27637} & \textbf{84\%} & \textbf{71\%} & \textbf{70\%} & \textbf{70\%} & \textbf{27718} & \textbf{13} & \textbf{0.29\%} & \textbf{27635} & \textbf{211}  & \textbf{22.4} & \textbf{28184} & \textbf{1.98\%}\\
\bottomrule
    \end{tabular}}
  \end{center}
  \caption{Results for instances with $N_c = 10$ and Small (S) intervals.}
  \label{tab:cvrp_results_1}
\end{table}

\newpage
\begin{table}[t]
  \begin{center}
    \resizebox{\textwidth}{!}{
    \begin{tabular}{ccrrrrrrrrrrrrr} 
      \toprule
      \multirow{2}{*}{$\mathcal{P}_{o}$} &
      \multirow{2}{*}{$\mathcal{P}_{m}$} & \multirow{2}{*}{\textbf{$\mathcal{S}_m$ cost}} &  \multirow{2}{*}{\textbf{$sim(\mathcal{S}_{o},\mathcal{S}_{m})$}} & \multicolumn{3}{c}{\textbf{ML model metrics}} & \multicolumn{3}{c}{\textbf{Edge-fixing}} & \multicolumn{3}{c}{\textbf{Exact BaP}} & \multicolumn{2}{c}{\textbf{DACT}}\\
      \cmidrule(rl){5-7} \cmidrule(rl){8-10} \cmidrule(rl){11-13} \cmidrule(rl){14-15}
       & & & & \textbf{TNR} & \textbf{TPR} & \textbf{Accuracy} & \textbf{Cost} & \textbf{Time (s)} & \textbf{Gap} & \textbf{Cost} & \textbf{Time (s)} & \textbf{Ratio} & \textbf{Cost} & \textbf{Gap}\\
      \midrule
      
      \multirow{5}{*}{\textbf{X-n101-k25}} & \textbf{\_10M\_1} & 27722 & 88\% & 83\% & 63\% & 73\% & 27969 & 36 & 0.89\% & 27722 & 544 & 15.1 & 28417 & 2.51\% \\ 
        & \textbf{\_10M\_2}  & 27636 & 76\% & 72\% & 66\% & 69\% & 27841 & 19 & 0.74\% & 27636 & 246 & 12.9 & 28156 & 1.88\% \\ 
        & \textbf{\_10M\_3} & 27411 & 88\% & 67\% & 61\% & 64\% & 27482 & 19 & 0.26\% & 27411 & 370 & 19.5 & 28143 & 2.67\% \\ 
        & \textbf{\_10M\_4} & 27655 & 69\% & 61\% & 65\% & 63\% & 27784 & 10 & 0.47\% & 27655 & 117 & 11.7 & 28286 & 2.28\% \\ 
        & \textbf{\_10M\_5} & 27605 & 96\% & 100\% & 59\% & 80\% & 27605 & 8 & 0.00\% & 27605 & 343 & 42.9 & 28489 & 3.20\% \\ 
        \textbf{Average} & - & \textbf{27606} & \textbf{83\%} & \textbf{77\%} & \textbf{63\%} & \textbf{70\%} & \textbf{27736} & \textbf{18} & \textbf{0.47\%} & \textbf{27606} & \textbf{324} & \textbf{20.4} & \textbf{28298} & \textbf{2.51\%} \\
\bottomrule
    \end{tabular}}
  \end{center}
  \caption{Results for instances with $N_c = 10$ and Medium (M) intervals.}
  \label{tab:cvrp_results_2}
\end{table}
\newpage
\begin{table}[t]
  \begin{center}
    \resizebox{\textwidth}{!}{
    \begin{tabular}{ccrrrrrrrrrrrrr} 
      \toprule
      \multirow{2}{*}{$\mathcal{P}_{o}$} &
      \multirow{2}{*}{$\mathcal{P}_{m}$} & \multirow{2}{*}{\textbf{$\mathcal{S}_m$ cost}} &  \multirow{2}{*}{\textbf{$sim(\mathcal{S}_{o},\mathcal{S}_{m})$}} & \multicolumn{3}{c}{\textbf{ML model metrics}} & \multicolumn{3}{c}{\textbf{Edge-fixing}} & \multicolumn{3}{c}{\textbf{Exact BaP}} & \multicolumn{2}{c}{\textbf{DACT}}\\
      \cmidrule(rl){5-7} \cmidrule(rl){8-10} \cmidrule(rl){11-13} \cmidrule(rl){14-15}
       & & & & \textbf{TNR} & \textbf{TPR} & \textbf{Accuracy} & \textbf{Cost} & \textbf{Time (s)} & \textbf{Gap} & \textbf{Cost} & \textbf{Time (s)} & \textbf{Ratio} & \textbf{Cost} & \textbf{Gap}\\
      \midrule
      
      \multirow{5}{*}{\textbf{X-n101-k25}} & \textbf{\_10L\_1} & 27781 & 83\% & 64\% & 70\% & 67\% & 28176 & 85 & 1.42\% & 27781 & 123 & 1.4 & 28420 & 2.30\% \\ 
        & \textbf{\_10L\_2} & 27493 & 61\% & 66\% & 86\% & 76\% & 27786 & 15 & 1.07\%  & 27493 & 154 & 10.3 & 27812 & 1.16\% \\ 
        & \textbf{\_10L\_3} & 27490 & 80\% & 78\% & 77\% & 78\% & 27526 & 15 & 0.13\%  & 27490 & 245 & 16.3 & 28297 & 2.94\% \\ 
        & \textbf{\_10L\_4} & 27775 & 64\% & 60\% & 78\% & 69\% & 28665 & 424 & 3.20\% & 27759 & 788 & 1.9 & 28248 & 1.70\% \\ 
        & \textbf{\_10L\_5} & 27717 & 64\% & 45\% & 74\% & 60\% & 28472 & 71 & 2.72\% & 27717 & 683 & 9.6 & 28210 & 1.78\% \\ 
        \textbf{Average} & - & \textbf{27651} & \textbf{70\%} & \textbf{63\%} & \textbf{77\%} & \textbf{70\%} & \textbf{28125} & \textbf{122} & \textbf{1.71\%} & \textbf{27648} & \textbf{399} & \textbf{7.9} & \textbf{28142} & \textbf{1.98\%}\\
\bottomrule
    \end{tabular}}
  \end{center}
  \caption{Results for instances with $N_c = 10$ and Large (L) intervals.}
  \label{tab:cvrp_results_3}
\end{table}
\newpage
\begin{table}[t]
  \begin{center}
    \resizebox{\textwidth}{!}{
    \begin{tabular}{ccrrrrrrrrrrrrr} 
      \toprule
      \multirow{2}{*}{$\mathcal{P}_{o}$} &
      \multirow{2}{*}{$\mathcal{P}_{m}$} & \multirow{2}{*}{\textbf{$\mathcal{S}_m$ cost}} &  \multirow{2}{*}{\textbf{$sim(\mathcal{S}_{o},\mathcal{S}_{m})$}} & \multicolumn{3}{c}{\textbf{ML model metrics}} & \multicolumn{3}{c}{\textbf{Edge-fixing}} & \multicolumn{3}{c}{\textbf{Exact BaP}} & \multicolumn{2}{c}{\textbf{DACT}}\\
      \cmidrule(rl){5-7} \cmidrule(rl){8-10} \cmidrule(rl){11-13} \cmidrule(rl){14-15}
       & & & & \textbf{TNR} & \textbf{TPR} & \textbf{Accuracy} & \textbf{Cost} & \textbf{Time (s)} & \textbf{Gap} & \textbf{Cost} & \textbf{Time (s)} & \textbf{Ratio} & \textbf{Cost} & \textbf{Gap}\\
      \midrule

      \multirow{5}{*}{\textbf{X-n101-k25}} & \textbf{\_20S\_1} & 27524  & 73\% & 75\% & 82\% & 79\% & 27761 & 13 & 0.86\% & 27524 & 140 & 10.8 & 28402 & 3.19\% \\ 
        & \textbf{\_20S\_2} & 27543 & 96\% & 75\% & 70\% & 73\% & 27543 & 27 & 0.00\% & 27543 & 328 & 12.1 & 28145 & 2.19\% \\ 
        & \textbf{\_20S\_3} & 27461 & 84\% & 87\% & 75\% & 81\% & 27551 & 26 & 0.33\% & 27461 & 227 & 8.7 & 28137 & 2.46\% \\ 
        & \textbf{\_20S\_4} & 27522 & 88\% & 64\% & 69\% & 67\% & 27709 & 44 & 0.68\% & 27522 & 288 & 6.5 & 27980 & 1.66\% \\ 
        & \textbf{\_20S\_5} & 27382 & 73\% & 62\% & 78\% & 70\% & 27575 & 36 & 0.70\% & 27382 & 171 & 4.8 & 28436 & 3.85\% \\ 
        \textbf{Average} & - & \textbf{27486} & \textbf{83\%} & \textbf{73\%} & \textbf{75\%} & \textbf{74\%} & \textbf{27628} & \textbf{29} & \textbf{0.51\%} & \textbf{27486} & \textbf{231} & \textbf{8.6} & \textbf{28220} & \textbf{2.67\%} \\
\bottomrule
    \end{tabular}}
  \end{center}
  \caption{Results for instances with $N_c = 20$ and Small (S) intervals.}
  \label{tab:cvrp_results_4}
\end{table}
\newpage
\begin{table}[t]
  \begin{center}
    \resizebox{\textwidth}{!}{
    \begin{tabular}{ccrrrrrrrrrrrrr} 
      \toprule
      \multirow{2}{*}{$\mathcal{P}_{o}$} &
      \multirow{2}{*}{$\mathcal{P}_{m}$} & \multirow{2}{*}{\textbf{$\mathcal{S}_m$ cost}} &  \multirow{2}{*}{\textbf{$sim(\mathcal{S}_{o},\mathcal{S}_{m})$}} & \multicolumn{3}{c}{\textbf{ML model metrics}} & \multicolumn{3}{c}{\textbf{Edge-fixing}} & \multicolumn{3}{c}{\textbf{Exact BaP}} & \multicolumn{2}{c}{\textbf{DACT}}\\
      \cmidrule(rl){5-7} \cmidrule(rl){8-10} \cmidrule(rl){11-13} \cmidrule(rl){14-15}
       & & & & \textbf{TNR} & \textbf{TPR} & \textbf{Accuracy} & \textbf{Cost} & \textbf{Time (s)} & \textbf{Gap} & \textbf{Cost} & \textbf{Time (s)} & \textbf{Ratio} & \textbf{Cost} & \textbf{Gap}\\
      \midrule
      
      \multirow{5}{*}{\textbf{X-n101-k25}} &  \textbf{\_20M\_1} & 27709 & 63\% & 77\% & 75\% & 76\% & 28105 & 47 & 1.43\% & 27709 & 550 & 11.7 & 28385 & 2.44\% \\ 
        & \textbf{\_20M\_2} & 27440 & 65\% & 68\% & 73\% & 71\% & 27646 & 64 & 0.75\% & 27440 & 615 & 9.6 & 27696 & 0.93\% \\ 
        & \textbf{\_20M\_3} & 27530 & 70\% & 81\% & 72\% & 77\% & 27559 & 24 & 0.11\% & 27530 & 983 & 41.0 & 28193 & 2.41\% \\ 
        & \textbf{\_20M\_4} & 27656 & 74\% & 70\% & 65\% & 68\% & 27767 & 11 & 0.40\% & 27656 & 217 & 19.7 & 28413 & 2.74\% \\ 
        & \textbf{\_20M\_5} & 27226 & 67\% & 76\% & 75\% & 76\% & 27393 & 15 & 0.61\% & 27226 & 223 & 14.9 & 27694 & 1.72\% \\ 
        \textbf{Average} & - & \textbf{27512} & \textbf{68\%} & \textbf{74\%} & \textbf{72\%} & \textbf{73\%} & \textbf{27694} & \textbf{32} & \textbf{0.66\%} & \textbf{27512} & \textbf{518} & \textbf{19.4} & \textbf{28076} & \textbf{2.05\%} \\
\bottomrule
    \end{tabular}}
  \end{center}
  \caption{Results for instances with $N_c = 20$ and Medium (M) intervals.}
  \label{tab:cvrp_results_5}
\end{table}

\newpage
\begin{table}[t]
  \begin{center}
    \resizebox{\textwidth}{!}{
    \begin{tabular}{ccrrrrrrrrrrrrr} 
      \toprule
      \multirow{2}{*}{$\mathcal{P}_{o}$} &
      \multirow{2}{*}{$\mathcal{P}_{m}$} & \multirow{2}{*}{\textbf{$\mathcal{S}_m$ cost}} &  \multirow{2}{*}{\textbf{$sim(\mathcal{S}_{o},\mathcal{S}_{m})$}} & \multicolumn{3}{c}{\textbf{ML model metrics}} & \multicolumn{3}{c}{\textbf{Edge-fixing}} & \multicolumn{3}{c}{\textbf{Exact BaP}} & \multicolumn{2}{c}{\textbf{DACT}}\\
      \cmidrule(rl){5-7} \cmidrule(rl){8-10} \cmidrule(rl){11-13} \cmidrule(rl){14-15}
       & & & & \textbf{TNR} & \textbf{TPR} & \textbf{Accuracy} & \textbf{Cost} & \textbf{Time (s)} & \textbf{Gap} & \textbf{Cost} & \textbf{Time (s)} & \textbf{Ratio} & \textbf{Cost} & \textbf{Gap}\\
      \midrule
      
      \multirow{5}{*}{\textbf{X-n101-k25}} & \textbf{\_20L\_1} & 27463 & 67\% & 70\% & 63\% & 67\% & 27637 & 107 & 0.63\% & 27463 & 374 & 3.5 & 28221 & 2.76\% \\ 
        & \textbf{\_20L\_2} & 27507 & 66\% & 69\% & 70\% & 70\% & 27641 & 93 & 0.49\% & 27507 & 196 & 2.1 & 28158 & 2.37\% \\ 
        & \textbf{\_20L\_3} & 27789 & 51\% & 76\% & 85\% & 81\% & 28245 & 50 & 1.64\% & 27789 & 100 & 2.0 & 28369 & 2.09\% \\ 
        & \textbf{\_20L\_4} & 27837 & 53\% & 65\% & 76\% & 71\% & 28097 & 7 & 0.93\% & 27837 & 707 & 101.0 & 28276 & 1.58\% \\ 
        & \textbf{\_20L\_5} & 27631 & 59\% & 70\% & 73\% & 72\% & 27733 & 21 & 0.37\% & 27631 & 955 & 45.5 & 28194 & 2.04\% \\ 
        \textbf{Average} & - & \textbf{27645} & \textbf{59\%} & \textbf{70\%} & \textbf{73\%} & \textbf{72\%} & \textbf{27871} & \textbf{56} & \textbf{0.81\%} & \textbf{27645} & \textbf{466} & \textbf{30.8} & \textbf{28244} & \textbf{2.17\%}\\
\bottomrule
      
    \end{tabular}}
  \end{center}
  \caption{Results for instances with $N_c = 20$ and Large (L) intervals.}
  \label{tab:cvrp_results_6}
\end{table}
\newpage
\begin{table}[t]
  \begin{center}
    \resizebox{\textwidth}{!}{
    \begin{tabular}{ccrrrrrrrrrrrrr} 
      \toprule
      \multirow{2}{*}{$\mathcal{P}_{o}$} &
      \multirow{2}{*}{$\mathcal{P}_{m}$} & \multirow{2}{*}{\textbf{$\mathcal{S}_m$ cost}} &  \multirow{2}{*}{\textbf{$sim(\mathcal{S}_{o},\mathcal{S}_{m})$}} & \multicolumn{3}{c}{\textbf{ML model metrics}} & \multicolumn{3}{c}{\textbf{Edge-fixing}} & \multicolumn{3}{c}{\textbf{Exact BaP}} & \multicolumn{2}{c}{\textbf{DACT}}\\
      \cmidrule(rl){5-7} \cmidrule(rl){8-10} \cmidrule(rl){11-13} \cmidrule(rl){14-15}
       & & & & \textbf{TNR} & \textbf{TPR} & \textbf{Accuracy} & \textbf{Cost} & \textbf{Time (s)} & \textbf{Gap} & \textbf{Cost} & \textbf{Time (s)} & \textbf{Ratio} & \textbf{Cost} & \textbf{Gap}\\
      \midrule
      
      \multirow{5}{*}{\textbf{X-n101-k25}} & \textbf{\_30S\_1} & 27645 & 73\% & 79\% & 69\% & 74\% & 27759 & 40 & 0.41\% & 27645 & 247 & 6.2 & 28109 & 1.68\% \\ 
        & \textbf{\_30S\_2} & 27709 & 71\% & 72\% & 72\% & 72\% & 27838 & 33 & 0.47\% & 27709 & 254 & 7.7 & 28579 & 3.14\% \\ 
        & \textbf{\_30S\_3} & 27485 & 78\% & 77\% & 66\% & 72\% & 27606 & 89 & 0.44\% & 27485 & 256 & 2.9 & 28104 & 2.25\% \\ 
        & \textbf{\_30S\_4} & 27620 & 84\% & 85\% & 66\% & 76\% & 27712 & 8 & 0.33\% & 27620 & 47 & 5.9 & 28422 & 2.90\% \\ 
        & \textbf{\_30S\_5} & 27350 & 68\% & 87\% & 79\% & 83\% & 27431 & 13 & 0.30\% & 27350 & 87 & 6.7 & 27883 & 1.95\% \\ 
        \textbf{Average} & - & \textbf{27562} & \textbf{75\%} & \textbf{80\%} & \textbf{70\%} & \textbf{75\%} & \textbf{27669} & \textbf{37} & \textbf{0.39\%} & \textbf{27562} & \textbf{178} & \textbf{5.9} & \textbf{28219} & \textbf{2.38\%}\\
\bottomrule

    \end{tabular}}
  \end{center}
  \caption{Results for instances with $N_c = 30$ and Small (S) intervals.}
  \label{tab:cvrp_results_7}
\end{table}

\newpage
\begin{table}[t]
  \begin{center}
    \resizebox{\textwidth}{!}{
    \begin{tabular}{ccrrrrrrrrrrrrr} 
      \toprule
      \multirow{2}{*}{$\mathcal{P}_{o}$} &
      \multirow{2}{*}{$\mathcal{P}_{m}$} & \multirow{2}{*}{\textbf{$\mathcal{S}_m$ cost}} &  \multirow{2}{*}{\textbf{$sim(\mathcal{S}_{o},\mathcal{S}_{m})$}} & \multicolumn{3}{c}{\textbf{ML model metrics}} & \multicolumn{3}{c}{\textbf{Edge-fixing}} & \multicolumn{3}{c}{\textbf{Exact BaP}} & \multicolumn{2}{c}{\textbf{DACT}}\\
      \cmidrule(rl){5-7} \cmidrule(rl){8-10} \cmidrule(rl){11-13} \cmidrule(rl){14-15}
       & & & & \textbf{TNR} & \textbf{TPR} & \textbf{Accuracy} & \textbf{Cost} & \textbf{Time (s)} & \textbf{Gap} & \textbf{Cost} & \textbf{Time (s)} & \textbf{Ratio} & \textbf{Cost} & \textbf{Gap}\\
      \midrule

      \multirow{5}{*}{\textbf{X-n101-k25}} &  \textbf{\_30M\_1} & 27849 & 57\% & 66\% & 75\% & 71\% & 28118 & 110 & 0.97\% & 27830 & 244 & 2.2 & 28118 & 0.97\% \\ 
        & \textbf{\_30M\_2} & 27660 & 60\% & 67\% & 72\% & 70\% & 27964 & 66 & 1.10\% & 27660 & 91 & 1.4 & 28338 & 2.45\% \\
        & \textbf{\_30M\_3} & 27699 & 69\% & 69\% & 72\% & 71\% & 27718 & 8 & 0.07\% & 27699 & 89 & 11.1 & 28401 & 2.53\% \\ 
        & \textbf{\_30M\_4} & 27724 & 60\% & 67\% & 71\% & 69\% & 27843 & 67 & 0.43\% & 27724 & 149 & 2.2 & 28388 & 2.40\% \\ 
        & \textbf{\_30M\_5} & 27340 & 68\% & 70\% & 76\% & 73\% & 27375 & 6 & 0.13\% & 27340 & 642 & 107.0 & 28139 & 2.92\% \\ 
        \textbf{Average} & - & \textbf{27654} & \textbf{63\%} & \textbf{68\%} & \textbf{73\%} & \textbf{71\%} & \textbf{27804} & \textbf{51} & \textbf{0.54\%} & \textbf{27651} & \textbf{243} & \textbf{24.8} & \textbf{28277} & \textbf{2.25\%}\\
\bottomrule
      
    \end{tabular}}
  \end{center}
  \caption{Results for instances with $N_c = 30$ and Medium (M) intervals.}
  \label{tab:cvrp_results_8}
\end{table}
\newpage

\begin{table}[t]
  \begin{center}
    \resizebox{\textwidth}{!}{
    \begin{tabular}{ccrrrrrrrrrrrrr} 
      \toprule
      \multirow{2}{*}{$\mathcal{P}_{o}$} &
      \multirow{2}{*}{$\mathcal{P}_{m}$} & \multirow{2}{*}{\textbf{$\mathcal{S}_m$ cost}} &  \multirow{2}{*}{\textbf{$sim(\mathcal{S}_{o},\mathcal{S}_{m})$}} & \multicolumn{3}{c}{\textbf{ML model metrics}} & \multicolumn{3}{c}{\textbf{Edge-fixing}} & \multicolumn{3}{c}{\textbf{Exact BaP}} & \multicolumn{2}{c}{\textbf{DACT}}\\
      \cmidrule(rl){5-7} \cmidrule(rl){8-10} \cmidrule(rl){11-13} \cmidrule(rl){14-15}
       & & & & \textbf{TNR} & \textbf{TPR} & \textbf{Accuracy} & \textbf{Cost} & \textbf{Time (s)} & \textbf{Gap} & \textbf{Cost} & \textbf{Time (s)} & \textbf{Ratio} & \textbf{Cost} & \textbf{Gap}\\
      \midrule
      
      \multirow{5}{*}{\textbf{X-n101-k25}} & \textbf{\_30L\_1} & 27736 & 61\% & 76\% & 75\% & 76\% & 28298 & 47 & 2.03\% & 27736 & 555 & 11.8 & 28486 & 2.70\% \\ 
        & \textbf{\_30L\_2} & 27465 & 72\% & 84\% & 71\% & 78\% & 27590 & 13 & 0.46\% & 27465 & 300 & 23.1 & 28045 & 2.11\% \\ 
        & \textbf{\_30L\_3} & 27809 & 66\% & 70\% & 70\% & 70\% & 28281 & 844 & 1.70\% & 27809 & 126 & 0.1 & 28438 & 2.26\% \\ 
        & \textbf{\_30L\_4} & 27809 & 65\% & 70\% & 69\% & 70\% & 28156 & 12 & 1.25\% & 27809 & 348 & 29.0 & 28358 & 1.97\% \\ 
        & \textbf{\_30L\_5} & 27267 & 65\% & 82\% & 73\% & 78\% & 27321 & 11 & 0.20\% & 27267 & 364 & 33.1 & 27709 & 1.62\% \\ 
        \textbf{Average} & - & \textbf{27617} & \textbf{66\%} & \textbf{76\%} & \textbf{72\%} & \textbf{74\%} & \textbf{27929} & \textbf{185} & \textbf{1.12\%} & \textbf{27617} & \textbf{339} & \textbf{19.4} & \textbf{28207} & \textbf{2.13\%}\\
\bottomrule
        
    \end{tabular}}
  \end{center}
  \caption{Results for instances with $N_c = 30$ and Large (L) intervals.}
  \label{tab:cvrp_results_9}
\end{table}

\end{document}